\documentclass[11pt]{article}
\usepackage{amsmath,amssymb}
\usepackage{CJK,indentfirst}
\usepackage{amssymb}
\begin{document}
\begin{CJK*}{GBK}{fs}
\baselineskip=20pt
\parskip=1mm
\parindent=20pt
\voffset=-1 true cm \setcounter{page} {1}
\numberwithin{equation}{section}
\newtheorem{Theorem}{\bf Theorem}[section]
\newtheorem{Prop}{\bf Proposition}[section]
\newtheorem{Lem}{\bf Lemma}[section]
\newtheorem{Def}{\bf Definition}[section]
\newtheorem{Coro}{\bf Corollary}[section]

\title{Existence and concentration of ground state solution to a 
 nonlocal Schr\"{o}dinger  
 equation  }

\author{ Anmin Mao$^a$  and Qian Zhang$^a$$^{b,}$\thanks{Corresponding author:\ zhangqianmath@163.com(Q.Zhang)}\\
\it \small$^a$ School of Mathematical Sciences, Qufu Normal University,  
Shandong  273165, P.R. China \\
\it \small$^b$ College of Mathematics and Informatics, 
Fujian Normal University, Fuzhou 350117, P.R. China\\ 
\small maoam@163.com(A.Mao) \ \ zhangqianmath@163.com(Q.Zhang)}
\date{}
\maketitle

{\bf \noindent Abstract}\quad 
We study a class of   Schr\"{o}dinger-Kirchhoff system involving critical exponent.  We aim to find suitable conditions to assure the existence of
     a positive ground state solution of Nehari-Poho\u{z}aev type $u_{\varepsilon}$ with exponential decay at infinity for $\varepsilon$ and $ u_{\varepsilon}$ concentrates around a global minimum point of $ V$ as $ \varepsilon\rightarrow0^{+}.$  The nonlinear term includes the nonlinearity $f(u)\sim|u|^{p-1}u$ for the well-studied case $ p\in[3,5)$, and the less-studied case $p\in(2,3)$.



\section{Introduction and main results}

 \noindent This paper is concerned with the following critical Schr\"{o}dinger-Kirchhoff type system
$$\left\{\aligned &-\left(\varepsilon^{2}a+\varepsilon b\int_{\mathbb{R}^{3}}|\nabla u|^{2}\right)\Delta u+V(x)u=A(x)|u|^{p-1}u+B(x)|u|^{4}u, &x \in \mathbb R^{3},\\
&u\in H^{1}(\mathbb{R}^{3}), \ \ u>0, &x \in \mathbb{R}^{3},\endaligned  \right. \eqno(1.1)$$
where $\varepsilon, a, b>0$, $ p \in  (2, 5)$. Problem ($1.1)$ can be used to  describe  several physical phenomena, see \cite{1,2,18,19,23,25,38} for more physical background. If we set  $\varepsilon=1,V(x)\equiv0,$ and $\mathbb R^{3}$ are replaced by bounded domain $\Omega,$ then problem $(1.1)$ are related to the following  Kirchhoff type problem
$$\left\{\aligned &-\left(a+b\int_{\mathbb{R}^{3}}|\nabla u|^{2}\right)\Delta u=f(x,u),&x \in \Omega,\\
&u=0,&x \in \partial \Omega,\endaligned \right.   $$
which is proposed by Kirchhoff in $\cite{1}$ and the solvability  has been well studied   only after Lions $\cite{2}$ introduced an abstract framework to this problem.
Using the minimax methods together with invariant sets of descent flow, Mao and Zhang $\cite{9}$, Perera and Zhang $\cite{10}$ proved the existence of  sign-changing solutions. When $f$ is critical growth,   $\cite{8,11}$ obtained the existence of positive solution  by variational methods.

Mathematically, ($1.1)$ is a nonlocal problem as the appearance of the nonlocal term $\int_{\mathbb{R}^{3}}|\nabla u|^{2}$  implies that ($1.1)$  is not a pointwise identity. This causes some mathematical difficulties which make the study of ($1.1)$  particularly interesting. Recently, there has been increasing attention to system
$$\left\{\aligned &-\left( a+  b\int_{\mathbb{R}^{3}}|\nabla u|^{2}\right)\Delta u+V(x)u=f(x,u), &x\in \mathbb R^{N},\\
&u\in H^{1}(\mathbb R^{N}),&x\in \mathbb R^{N},\endaligned\right. \eqno(1.2)$$
 $V$ is always assumed to satisfy
$$(V)  \ \ \  \ \ \ \ \inf\limits_{x\in \mathbb R^{N}}V(x)\geq a>0  \  and \ for \ each  \ M>0,\  meas\{x\in \mathbb R^{N}:V(x)\leq M\}<+\infty.\ \ \ \ \ \ \ \ \ \ \ \  $$
Such kind of hypothesis has been widely used and was firstly introduced by Bartsch and Wang $\cite{18}$ to ensure the compactness of embedding of $E\triangleq \{u\in H^{1}(\mathbb R^{3})\ |\ \int_{\mathbb R^{3}}V(x)|u|^{2}<\infty\}\hookrightarrow L^{s}(\mathbb R^{3}),$ where $2\leq s<6,$ for example \cite{19,20}. Another important hypothesis on nonlinear term is the following monotonicity condition:
$$ (F)  \ \ \ \ \ \  \ \ \ \  \frac{f(u)}{u^{3}} \   is \ strictly \ increasing \ for  \    u>0,\ \ \ \ \ \ \ \ \ \ \ \ \ \ \ \ \ \ \ \ \ \ \ \ \ \ \ \ \ \ \ \ \ \ \ \ \ \ \ \ \ \ \ \ \ \ \ \ \ \ \ \ \ \ \ \ \ \ \ \ \ \  \ \  \ \  \ \  \ \ \ \ \ \  \ \ $$
which is key to prove the compactness of $(PS)$ sequence. In the subcritical case, by using the variational methods and some estimates, He and Zou $\cite{21}$ studied the concentration of positive ground state solution of Nehari type of (1.2) when  $f$ satisfies the $(F)$. When $f(x,u)$  in (1.2) involves critical exponent, one has to face lots of difficulties caused by  the critical  exponent, in such case, $(F)$ also plays an important role in  variational  approach, see \cite{22,23,24,43,46}.

 Regarding the existence of positive ground state solutions for critical Kirchhoff type problem, to the best of our knowledge, there are very few results in the context.  By variational methods,  Wang et al. \cite{22} and Li and Ye \cite{31}  proved the existence of positive ground state solutions to subcritical problems, it is worthy pointing out that  conditions $(V)$ and $(F)$ are critical to their approach.

Inspired by above works, in this paper, we restrict our attention to a class of more general problem  $(1.1)$ involving critical exponent and study the existence and concentration of positive ground state solutions without coercive condition $(V)$ and monotonicity condition $(F)$. Moreover,  we also investigate the exponent decay of positive ground state solutions.

Before state our main results, we  make the some assumptions.\\
$(V_{1})$ $V\in C(\mathbb R^{3},\mathbb R^{+}).$\\
$(V_{2})$ For almost every $x\in \mathbb R^{3},$
$$0<\inf_{x\in \mathbb R^{3}}V(x)=:V_{0}\leq V(x)\leq\liminf_{|x|\rightarrow\infty} V(x)=:V_{\infty}<+\infty \ \hbox{and} \ V(x)\not\equiv V_{\infty}.$$
$(V_{3})$ $V$ is weakly differentiable and satisfies $(\nabla V(x),x)\in L^{\infty}(\mathbb R^{3})\cup L^{\frac{3}{2}}(\mathbb R^{3})$ and\\
\indent \indent \indent \indent \indent $V(x)-(\nabla V(x),x)\geq0$,\ \ $x\in \mathbb R^{3}.$\\
$(V_{4})$ There exists a positive constant $l<a$ such that
$$|(\nabla V(x),x)|\leq \frac{l}{2|x|^{2}},\ \ \ x\in \mathbb R^{3}\backslash\{0\}.$$
$(F_{1})$ $A$ is weakly differentiable and satisfies $(\nabla A(x),x)\in L^{\infty}(\mathbb R^{3})\cup L^{2}(\mathbb R^{3})$ and\\
\indent \indent \indent \indent \indent $(\nabla A(x),x)\geq0$,\ \ $x\in \mathbb R^{3}.$\\
$(F_{2})$ $B$ is weakly differentiable and satisfies $(\nabla B(x),x)\in L^{\infty}(\mathbb R^{3})\cup L^{2}(\mathbb R^{3})$ and\\
\indent \indent \indent \indent \indent $(\nabla B(x),x)\geq0$,\ \ $x\in \mathbb R^{3}.$\\
$(F_{3})$  $A\in C(\mathbb R^{3},\mathbb R),\lim\limits_{|x|\rightarrow\infty}A(x)=A_{\infty}\in(0,\infty) $ and $A(x)>A_{\infty}>0$ for $x\in \mathbb R^{3}.$\\
$(F_{4})$ $B\in C(\mathbb R^{3},\mathbb R),\lim\limits_{|x|\rightarrow\infty}B(x)=B_{\infty}\in(0,\infty)$ and $B(x)>B_{\infty}>0$ for $x\in \mathbb R^{3}.$\\
$(F_{5})$ $\Lambda \cap \Lambda_{1}\cap \Lambda_{2}\neq\emptyset,$ where
$$\Lambda\triangleq\{x\in \mathbb R^{3}:V(x)=V_{0}\},\ \ \ \ \ \ \ \ \ \  \ \ \ \ \ \ \  \   $$
$$\Lambda_{1}\triangleq\{x\in \mathbb R^{3}:A(x)=A_{0}\triangleq\max_{x\in \mathbb R^{3}} A(x)\},$$
$$\Lambda_{2}\triangleq\{x\in \mathbb R^{3}: B(x)=B_{0}\triangleq\max_{x\in \mathbb R^{3}}B(x)\}.$$

\vskip4pt

Note that conditions $(V_{3})$, $(F_{1})$ and $(F_{2})$ are usually known as the Poho\u{z}aev type conditions. About their applications, we refer the reader to $\cite{25,26,31}$. $(F_{3})$-$(F_{5})$ once appeared in the study of nonlinear Schr\"{o}dinger and Schr\"{o}dinger-Possion system in $\cite{28,29,30}$.
\vskip4pt

Next, we give some notations. Throughout this paper, let $H^{1}(\mathbb R^{3})$ be the usual Sobolev space equipped with the norm $\|u\|_{H}:=\left(\int_{\mathbb R^{3}}a|\nabla u|^{2}+|u|^{2}\right)^{\frac{1}{2}},$ $(\int_{\mathbb R^{3}}|\nabla u|^{2})^{\frac{1}{2}}$ the norm of $D^{1,2}(\mathbb R^{3}),$ and $|\cdot|_{q}$ the usual Lebesgue space $L^{q}(\mathbb R^{3})$ norm.

By  the change of variable $\varepsilon z=x,$ we   rewrite $(1.1)$ as
$$-\left(a+b\int_{\mathbb{R}^{3}}|\nabla u|^{2}\right)\Delta u+V(\varepsilon x)u=A(\varepsilon x)|u|^{p-1}u+B(\varepsilon x)|u|^{4}u,\ \ x\in\mathbb R^{3}.\eqno(K_{\varepsilon})$$
For $\rho>0$ and $z^{*}\in \mathbb R^{3},\ B_{\rho}(z^{*}):=\{z\in \mathbb R^{3}:|z^{*}-z|\leq\rho\}.$ For any $\varepsilon >0,$  let $ E_{\varepsilon}=\{u\in H^{1}(\mathbb R^{3}):\int_{\mathbb R^{3}}V(\varepsilon x)|u|^{2}<\infty\}$ denote the Hilbert space endowed with inner product and norm $\|u\|_{\varepsilon}:=\langle u,u\rangle_{\varepsilon}^{\frac{1}{2}}=\left(\int_{\mathbb R^{3}}a|\nabla u|^{2}+V(\varepsilon x)|u|^{2}\right)^{\frac{1}{2}}.$ For simplicity, $ C, C_{i}$ will denote various positive constant.

$(K_{\varepsilon})$ is called variational which means that its weak solutions are the critical points of the functional $I_{\varepsilon}:u\in E_{\varepsilon}\rightarrow \mathbb R,$
$$ I_{\varepsilon}(u):=\frac{1}{2}\|u\|_{\varepsilon}^{2}+\frac{b}{4}\left(\int_{\mathbb R^{3}}|\nabla u|^{2}\right)^{2}-\frac{1}{p+1}\int_{\mathbb R^{3}}A(\varepsilon x)|u|^{p+1}-\frac{1}{6}\int_{\mathbb R^{3}}B(\varepsilon x)|u|^{6},\eqno(1.3)$$
which is a well-defined $C^{1}$ functional (see \cite{24}) and for all $\varphi\in E_{\varepsilon},$
$$ \aligned \langle I_{\varepsilon}'(u),\varphi\rangle&=\int_{\mathbb R^{3}}(a\nabla u\nabla \varphi+V(\varepsilon x)u\varphi)+b\int_{\mathbb R^{3}}|\nabla u|^{2}\int_{\mathbb R^{3}}\nabla u\nabla \varphi-\int_{\mathbb R^{3}}A(\varepsilon x)|u|^{p-1}u\varphi\\
&\ \ \ -\int_{\mathbb R^{3}}B(\varepsilon x)|u|^{4}u\varphi.\endaligned \eqno(1.4)$$

 Our main result  reads as follows.
\vskip4pt

{\it\noindent {\bf Theorem 1.1.} Under the assumptions $(V_{1})$-$(V_{4})$, $(F_{1})$-$(F_{5}),$ if $p\in(2,5)$, then \\
(1) there exists $ \varepsilon^{*}>0$ such that problem $(1.1)$ has a positive ground state solution of Nehari-Poho\u{z}aev type $ u_{\varepsilon}$ for $ \varepsilon\in(0,\varepsilon^{*})$;\\
(2) $ u_{\varepsilon}$ possesses a global maximum point $x_{\varepsilon}$ in $ \mathbb R^{3}$ such that $ \lim\limits_{\varepsilon\rightarrow0^{+}}V(x_{\varepsilon})=V_{0}.$
And, there are positive constants $ C,c$ such that
$ u_{\varepsilon}(x)\leq C\exp\left(-c\frac{|x-x_{\varepsilon}|}{\varepsilon} \right) \hbox{ \ for\ all } x\in \mathbb R^{3}. $}

\noindent {\bf Remark 1.2.}  Due to the emergence of critical terms  and the lack of coercive condition $(V)$ and monotonicity condition $(F)$ once used in \cite{21,22,23,24,43,46},  the approaches in above works can not be applied here.

\vskip4pt

\noindent {\bf Remark 1.3.} Different from  \cite{21,22,23,24,43,46},  we not only prove the existence of singed ground state solution but also figure out its exponent decay.  Hence, our results improve  some previous works.
\vskip4pt

\noindent {\bf Remark 1.4.} One of our main purpose is  to  find an unified approach for the full subcritical range of $p\in (2,5)$,  the approaches   by using Nehari manifold in \cite{21,22,23,24} and \cite{43,46} can not be applied here because that when applied to the monomial nonlinearity $f(u)=|u|^{p-1}u$,  the approaches in above mentioned papers are only valid for $p\in [3,5)$ and it is difficult to prove the boundedness of the minimizing sequence in our case.  Here, motivated by \cite{31}, we adopt  constrained minimization on a manifold and take the minimum on the new manifold.   Such a new manifold is obtained by combining the  Poho\u{z}aev identity and the Nehari manifold.
\vskip4pt

To complete the proof of theorem 1.1, we use  an indirect approach by using  Jeanjean's monotonicity trick in $\cite{27}$.
\vskip4pt

\noindent {\bf Lemma 1.5.} (\cite{27}) {\it Let $ (X,\|\cdot\|)$ be a Banach space and $ J\in \mathbb R^{+}$ be an interval. Consider a family of $C^{1}$ functionals on $X$ of the form
$$ \Phi_{\sigma}(u)=C(u)-\sigma D(u),\ \hbox{ for all } \ \sigma\in J,$$
with $D(u)\geq0$ and either $ C(u)\rightarrow +\infty$ or $ D(u)\rightarrow +\infty,$ as $\|u\|\rightarrow \infty.$ Assume that there are two points $ v_{1}, v_{2}\in X$ such that
$$c_{\sigma}=\inf_{\gamma\in \Gamma}\max_{s\in [0,1]}\Phi_{\sigma}(\gamma(s))>\max\{\Phi_{\sigma}(v_{1}), \Phi_{\sigma}(v_{2})\}, \ \ \hbox{ for any } \ \sigma\in J,$$
where
$$ \Gamma=\{\gamma\in C([0,1],X)|\ \gamma(0)=v_{1},\ \gamma(1)=v_{2}\}.$$
Then, for almost every $\sigma\in J,$ there is a bounded $(PS)_{c_{\sigma}}$ sequences in $ X.$}

Set $J=[\delta,1]$, where $ \delta\in(0,1)$. We consider a family of $C^{1}$ functionals $I_{\varepsilon,\sigma}(u),$ i.e., for any $\sigma\in[\delta,1]$,
$$ I_{\varepsilon,\sigma}(u)=\frac{1}{2}\|u\|_{\varepsilon}^{2}+\frac{b}{4}\left(\int_{\mathbb R^{3}}|\nabla u|^{2}\right)^{2}-\frac{\sigma}{p+1}\int_{\mathbb R^{3}}A(\varepsilon x)|u|^{p+1}-\frac{\sigma}{6}\int_{\mathbb R^{3}}B(\varepsilon x)|u|^{6}.$$
By $(V_{2}),(F_{3})$-$(F_{5})$ and Lemma 1.5, for a.e. $\sigma\in[\delta,1]$, there exists a bounded $(PS)_{c_{\sigma}}$ sequence in $H^{1}(\mathbb R^{3})$, denoted by $\{u_{n}\}$, where $c_{\sigma}$ is given below (see Lemma 3.4). The nonlinearity include $|u|^{p-1}u$ with $p\in(2,5)$ implies that the monotonicity of $\frac{|u|^{p-1}u}{u^{3}}$ does not always hold, so  we have to consider the following functional
$$ J_{\varepsilon,\sigma}(u)=\frac{a+b\varrho^{2}}{2}\int_{\mathbb R^{3}}|\nabla u|^{2}+\frac{1}{2}\int_{\mathbb R^{3}}V(\varepsilon x)|u|^{2}-\frac{\sigma}{p+1}\int_{\mathbb R^{3}}A(\varepsilon x)|u|^{p+1}-\frac{\sigma}{6}\int_{\mathbb R^{3}}B(\varepsilon x)|u|^{6},$$
  Hypothesis $(V_{2}),$ $(F_{3}),$ $(F_{4})$ can help us establish a global compactness lemma of critical type  (see Lemma 3.7 below) related to the functional $J_{\varepsilon,\sigma}$ and its limit  functional
$$ J_{\sigma}^{\infty}(u)=\frac{a+b\varrho^{2}}{2}\int_{\mathbb R^{3}}|\nabla u|^{2}+\frac{1}{2}\int_{\mathbb R^{3}}V_{\infty}|u|^{2}-\frac{\sigma}{p+1}\int_{\mathbb R^{3}}A_{\infty}|u|^{p+1}-\frac{\sigma}{6}\int_{\mathbb R^{3}}B_{\infty}|u|^{6}. $$
In what follows,  we shall consider  the  limit  problem of $(K_{\varepsilon})$, i.e.,
$$\left\{\aligned &-\left(a+b\int_{\mathbb R^{3}}|\nabla u|^{2}\right)\Delta u+V_{\infty}u=\sigma A_{\infty}|u|^{p-1}u+\sigma B_{\infty}|u|^{4}u, &x\in\mathbb R^{3}, \\
&u\in H^{1}(\mathbb R^{3}),\ u>0, &x\in\mathbb R^{3},\endaligned\right. \eqno(1.5)$$
and establish the existence of  the   least energy of the corresponding energy functional
$$ I_{\sigma}^{\infty}(u)=\frac{1}{2}\int_{\mathbb R^{3}}\left(a|\nabla u|^{2}+V_{\infty}|u|^{2}\right)+\frac{b}{4}\left(\int_{\mathbb R^{3}}|\nabla u|^{2}\right)^{2}-\frac{\sigma}{p+1}\int_{\mathbb R^{3}}A_{\infty}|u|^{p+1}-\frac{\sigma}{6}\int_{\mathbb R^{3}}B_{\infty}|u|^{6}.$$
To get the desired result, we need to prove the following result.

{\it\noindent {\bf Theorem 1.6.} (1.5) has a positive ground state solution of Nehari-Poho\u{z}aev type in $H^{1}(\mathbb R^{3})$ for every $2<p<5.$ }

As for problem (1.5), $I_{\sigma}^{\infty}$ does not always satisfy $(PS)_{c_{\sigma}}$ condition and it is difficult to get a ground state solution of Nehari type since $2<p<5$. For simplicity, we may assume that $V_{\infty}=A_{\infty}=B_{\infty}=\sigma\equiv1$ in (1.5), i.e.,
$$\left\{\aligned &-\left(a+b\int_{\mathbb R^{3}}|\nabla u|^{2}\right)\Delta u+u=|u|^{p-1}u+|u|^{4}u, & x\in\mathbb R^{3},\\
&u\in H^{1}(\mathbb R^{3}),\ u>0,  &x\in\mathbb R^{3},\endaligned\right. \eqno(1.6)$$
and the corresponding energy functional is
$$ I(u)=\frac{1}{2}\int_{\mathbb R^{3}}(a|\nabla u|^{2}+|u|^{2})+\frac{b}{4}\left(\int_{\mathbb R^{3}}|\nabla u|^{2}\right)^{2}-\frac{1}{p+1}\int_{\mathbb R^{3}}|u|^{p+1}-\frac{1}{6}\int_{\mathbb R^{3}}|u|^{6}.  \eqno(1.7)$$
To prove problem (1.6) has a positive ground state solution of Nehari type in $H^{1}(\mathbb R^{3})$ for every $p\in(2,5),$ as we describe before, the usual Nehari manifold is not suitable because it is difficult to prove the boundedness of the minimizing sequence. Inspired by [31], we   take the minimum on a new manifold which is constructed by combining the Nehari manifold and the corresponding Poho\u{z}aev type identity, which will be proved in Section 2 (see Lemma 2.1).  Finally, it follows from Theorem 1.6 and global compactness lemma of critical type that $(PS)_{c_{\sigma}}$ condition holds. Choosing a sequence $\{\sigma_{n}\}\subset[\delta,1]$ with $\sigma_{n}\rightarrow1$, there exist a sequence of nontrivial weak solutions $\{u_{\sigma_{n}}\}\subset H^{1}(\mathbb R^{3})$. We shall prove  $\{u_{\sigma_{n}}\}$ is a bounded $(PS)_{c_{1}}$ sequence for $I_{\varepsilon}=I_{\varepsilon,1}$ by using the Poho\u{z}aev identity, which yields the existence of positive ground state solution of Nehari type. Such  ground state solution  is also called of Nehari-Poho\u{z}aev type.

 The idea of the proof of Theorem 1.1 (2) is to use the Poho\u{z}aev type identity and some estimates.
Assume $0\in\Lambda\cap \Lambda_{1}\cap \Lambda_{2}$  and consider the associated autonomous problem
$$\left\{\aligned &-\left(a+b\int_{\mathbb{R}^{3}}|\nabla u|^{2}\right)\Delta u+V_{0}u=A_{0}|u|^{p-1}u+B_{0}|u|^{4}u, &x \in \mathbb{R}^{3},\\
&u\in H^{1}(\mathbb{R}^{3}),\ \  u>0, &x \in \mathbb{R}^{3},\endaligned\right. \eqno(K_{0})$$
 its corresponding functional
$$ I_{0}(u)=\frac{a}{2}\int_{\mathbb R^{3}}|\nabla u|^{2}+\frac{1}{2}\int_{\mathbb R^{3}}V_{0}|u|^{2}+\frac{b}{4}\left(\int_{\mathbb R^{3}}|\nabla u|^{2}\right)^{2}-\frac{1}{p+1}\int_{\mathbb R^{3}}A_{0}|u|^{p+1}-\frac{1}{6}\int_{\mathbb R^{3}}B_{0}|u|^{6},$$
  will be used to prove the concentration of positive ground state solutions. Motivated  but different from the work in \cite{24}, we do not need to prove the compactness of the minimizing sequence by using concentration-compactness principle,   thanks to the application of a new manifold, our proof becomes simpler.

\vskip4pt

 In Section 2, we state the variational framework of our problem and some preliminary results. Section 3 is devoted to the existence of positive ground state solution of Nehari-Poho\u{z}aev type.  And we investigate the concentration of ground state solution of Nehari-Poho\u{z}aev type and complete the proof of Theorem 1.1 in Section 4.

\section{Preliminaries and functional setting}

We begin with a Poho\u{z}aev identity with respect to problem $(K_{\varepsilon}),$ which plays a significant role in constructing a new manifold.

{\it\noindent {\bf Lemma 2.1.} Under the assumptions $(V_{1})$-$(V_{4})$ and $(F_{1})$-$(F_{5})$, if $u\in E_{\varepsilon}$ be a weak solution of problem $(K_{\varepsilon})$ and $p\in (2,5)$, then $u$ satisfies:
$$ \aligned 0=P_{\varepsilon}(u):=&\ \frac{a}{2}\int_{\mathbb R^{3}}|\nabla u|^{2}+\frac{3}{2}\int_{\mathbb R^{3}}V(\varepsilon x)|u|^{2}+\frac{\varepsilon}{2}\int_{\mathbb R^{3}}(\nabla V(\varepsilon x),x)|u|^{2}+\frac{b}{2}\left(\int_{\mathbb R^{3}}|\nabla u|^{2}\right)^{2}\\
&-\frac{3}{p+1}\int_{\mathbb R^{3}}A(\varepsilon x)|u|^{p+1}-\frac{\varepsilon}{p+1}\int_{\mathbb R^{3}}(\nabla A(\varepsilon x),x)|u|^{p+1}-\frac{1}{2}\int_{\mathbb R^{3}}B(\varepsilon x)|u|^{6}\\
&-\frac{\varepsilon}{6}\int_{\mathbb R^{3}}(\nabla B(\varepsilon x),x)|u|^{6}. \endaligned \eqno(2.1)$$}
\noindent {\it Proof.} It's easy to get the Poho\u{z}aev identity on a ball $B_{R}=\{x\in \mathbb R^{3}:\ |x|<R\}$:
$$\aligned&\frac{a}{2}\int_{B_{R}}|\nabla u|^{2}dx+\frac{b}{2}\left(\int_{B_{R}}|\nabla u|^{2}\right)^{2}dx+\frac{\varepsilon}{2}\int_{B_{R}}(\nabla V(\varepsilon x), x)|u|^{2}dx+\frac{3}{2}\int_{B_{R}}V(\varepsilon x)|u|^{2}dx\\
&\ -\frac{\varepsilon}{p+1}\int_{B_{R}}\left(\nabla A(\varepsilon x),x\right)|u|^{p+1}dx-\frac{3}{p+1}\int_{B_{R}}A(\varepsilon x)|u|^{p+1}dx-\frac{\varepsilon}{6}\int_{B_{R}}\left(\nabla B(\varepsilon x),x\right)|u|^{6}dx\\
&\ -\frac{1}{2}\int_{B_{R}}B(\varepsilon x)|u|^{6}dx\\
=&\ R\int_{\partial B_{R}}\bigg(-\left(a+b\int_{\mathbb R^{3}}|\nabla u|^{2}\right)\left(\frac{\partial u}{\partial{\bf n}}\right)^{2}+\frac{1}{2}\left(a+b\int_{\mathbb R^{3}}|\nabla u|^{2}\right)|\nabla u|^{2}+\frac{1}{2}V(\varepsilon x)|u|^{2}\\
&\ -\frac{1}{p+1}A(\varepsilon x)|u|^{p+1}-\frac{1}{6}B(\varepsilon x)|u|^{6}\bigg)\ dS. \endaligned \eqno(2.2)$$

Next we show that the right hand side in (2.2) converges to 0 for at least one suitably chosen sequence $R_{n}\rightarrow+\infty$. Since
$$ \aligned +\infty&>\int_{\mathbb R^{3}}\left|-\left(\frac{\partial u}{\partial{\bf n}}\right)^{2}+|\nabla u|^{2}+|u|^{2}-|u|^{p+1}-|u|^{6}\right|dx\\
&=\int_{0}^{+\infty}\left(\int_{\partial B_{R}}\left|-\left(\frac{\partial u}{\partial{\bf n}}\right)^{2}+|\nabla u|^{2}+|u|^{2}-|u|^{p+1}-|u|^{6}\right|\ dS\right)dR,\endaligned \eqno(2.3)$$
there exists a sequence $R_{n}\rightarrow+\infty$ such that
$$ R_{n}\int_{\partial B_{R_{n}}}\left|-\left(\frac{\partial u}{\partial{\bf n}}\right)^{2}+|\nabla u|^{2}+|u|^{2}-|u|^{p+1}-|u|^{6}\right|dS\rightarrow0 \ \ \hbox{as}\  n\rightarrow+\infty.$$
Taking limits as $n\rightarrow+\infty$ in (2.2) yields,
$$ \aligned &\frac{a}{2}\int_{\mathbb R^{3}}|\nabla u|^{2}+\frac{3}{2}\int_{\mathbb R^{3}}V(\varepsilon x)|u|^{2}+\frac{\varepsilon}{2}\int_{\mathbb R^{3}}(\nabla V(\varepsilon x),x)|u|^{2}+\frac{b}{2}\left(\int_{\mathbb R^{3}}|\nabla u|^{2}\right)^{2}\\
&-\frac{3}{p+1}\int_{\mathbb R^{3}}A(\varepsilon x)|u|^{p+1}-\frac{\varepsilon}{p+1}\int_{\mathbb R^{3}}(\nabla A(\varepsilon x),x)|u|^{p+1}-\frac{1}{2}\int_{\mathbb R^{3}}B(\varepsilon x)|u|^{6}\\
&-\frac{\varepsilon}{6}\int_{\mathbb R^{3}}(\nabla B(\varepsilon x),x)|u|^{6}=0. \endaligned $$ $\hfill\square$

In particular, if $V\equiv A\equiv B\equiv1$, we have
$$ P(u):=\frac{1}{2}a\int_{\mathbb R^{3}}|\nabla u|^{2}+\frac{3}{2}\int_{\mathbb R^{3}}|u|^{2}+\frac{1}{2}b\left(\int_{\mathbb R^{3}}|\nabla u|^{2}\right)^{2}-\frac{3}{p+1}\int_{\mathbb R^{3}}|u|^{p+1}-\frac{1}{2}\int_{\mathbb R^{3}}|u|^{6}=0.\eqno(2.4)$$

{\it\noindent {\bf Lemma 2.2.} If $p\in(2,5),$ then $I$ is not bounded from below.}\\
\noindent {\it Proof.} Let $u_{t}(x):=tu(t^{-1}x), t>0.$ Since $2<p<5,$ we have
$$ I(u_{t})=\frac{t^{3}}{2}a\int_{\mathbb R^{3}}|\nabla u|^{2}+\frac{t^{5}}{2}\int_{\mathbb R^{3}}|u|^{2}+\frac{t^{6}}{4}b\left(\int_{\mathbb R^{3}}|\nabla u|^{2}\right)^{2}-\frac{ t^{p+4}}{p+1}\int_{\mathbb R^{3}}|u|^{p+1}-\frac{t^{9}}{6}\int_{\mathbb R^{3}}|u|^{6}\rightarrow-\infty$$
as $ t\rightarrow+\infty$ for all $u\in H^{1}(\mathbb R^{3})\backslash\{0\}.$ $\hfill\square$

Lemma 2.2 shows that $ I$ has a mountain pass geometry around $0\in H^{1}(\mathbb R^{3})$ which is a local minimum of $I$. However, it is difficult for us to obtain the boundedness of the $(PS)$ sequence for $p\in(2,5).$ Therefore, we need to consider the constrained minimization on a suitable manifold.

 To give the definition of such a manifold, we need the following lemma.
\vskip4pt
{\it\noindent {\bf Lemma 2.3.} Let $ a_{i}\ (i=1,2,3,4,5)$ be positive constants and $p\in(2,5).$ If $ h(t)=a_{1}t^{3}+a_{2}t^{5}+a_{3}t^{6}-a_{4}t^{p+4}-a_{5}t^{9}$ for $ t\geq0.$ Then $h$ has a unique critical point which corresponds to its maximum.}\\
\noindent {\it Proof.} The technique is the same as in \cite[Lemma 3.3]{32} and we omit the proof here.   $\hfill\square$

If $u\in H^{1}(\mathbb R^{3})$ is a nontrivial critical point of $ I$, set
$$ \chi(t)\triangleq I(u_{t})=\frac{t^{3}}{2}a\int_{\mathbb R^{3}}|\nabla u|^{2}+\frac{t^{5}}{2}\int_{\mathbb R^{3}}|u|^{2}+\frac{t^{6}}{4}b\left(\int_{\mathbb R^{3}}|\nabla u|^{2}\right)^{2}-\frac{t^{p+4}}{p+1}\int_{\mathbb R^{3}}|u|^{p+1}-\frac{t^{9}}{6}\int_{\mathbb R^{3}}|u|^{6}.$$
Using Lemma 2.3 we know that $\chi$ has a unique critical point $ t_{0}>0$ corresponding to its maximum. Since $u$ is a solution to (1.6), we see that $t_{0}=1$ and $ \chi'(1)=0,$ i.e.,
$$ \aligned G(u)&:=\frac{3}{2}a\int_{\mathbb R^{3}}|\nabla u|^{2}+\frac{5}{2}\int_{\mathbb R^{3}}|u|^{2}+\frac{3}{2}b\left(\int_{\mathbb R^{3}}|\nabla u|^{2}\right)^{2}-\frac{p+4}{p+1}\int_{\mathbb R^{3}}|u|^{p+1}-\frac{3}{2}\int_{\mathbb R^{3}}|u|^{6}=0. \endaligned \eqno(2.5)$$
It is easy to see that $G(u)=\langle I'(u),u \rangle+P(u)$, where $P(u)$ is given in (2.4). Then we define
$$ M=\{u\in H^{1}(\mathbb R^{3})\backslash\{0\} : G(u)=0\}.$$

{\it\noindent {\bf Lemma 2.4.} Assume $ p\in(2,5).$ Then for any $0\neq u\in H^{1}(\mathbb R^{3}),$ there exists a unique $\hat{t}>0,$ such that $u_{\hat{t}}\in M,$ where $u_{\hat{t}}(x)=\hat{t}u(\hat{t}^{-1}x).$ Moreover, $I(u_{\hat{t}})=\max\limits_{t>0}I(u_{t}).$ }

\noindent {\it Proof.} The proof is standard and we omit it here. $\hfill\square$
\vskip3pt

{\it\noindent {\bf Lemma 2.5.} If $p\in(2,5)$, then $M$ is a natural $C^{1}$ manifold and every critical point of $ I|_{M}$ is a critical point of $I$ in $ H^{1}(\mathbb R^{3}).$ }

\noindent {\it Proof.} By Lemma 2.4, it is then easy to check that $M\neq\emptyset.$ The proof consists of four steps.

{\bf Step 1.} $0\not \in\partial M$.

Note that, for any $u\in M,$ using the Sobolev embedding inequality, choosing a number $\rho>0$, then there exist $ r>0,C_{1}, C_{2}>0$ such that
$$ \aligned G(u)&=\frac{3}{2}a\int_{\mathbb R^{3}}|\nabla u|^{2}+\frac{5}{2}\int_{\mathbb R^{3}}|u|^{2}+\frac{3}{2}b\left(\int_{\mathbb R^{3}}|\nabla u|^{2}\right)^{2}-\frac{p+4}{p+1}\int_{\mathbb R^{3}}|u|^{p+1}-\frac{3}{2}\int_{\mathbb R^{3}}|u|^{6}\\
&\geq\frac{3}{2}\|u\|_{H}^{2}-C_{1}\frac{p+4}{p+1}\|u\|_{H}^{p+1}-\frac{3}{2}C_{2}\|u\|_{H}^{6}>r>0,\endaligned$$
for $\|u\|_{H}=\rho $ small enough, so that $M, \partial M\subset H^{1}(\mathbb R^{3})\backslash B_{\rho}(0)$.

{\bf Step 2.} $\inf\limits_{M} I(u)>0$.

If $u\in M$, we have
$$  (p+4)I(u) =(p+4)I(u)-G(u)\geq\frac{p-1}{2}\|u\|_{H}^{2}>\frac{p-1}{2}\rho^{2},  \eqno(2.6)$$
which implies $\inf\limits_{M} I(u)>0$.

{\bf Step 3.} $M$ is a $C^{1}$ manifold.

Since $G(u)$ is a $C^{1}$ functional, in order to prove $M$ is a $C^{1}$ manifold, it suffices to prove that $G'(u)\neq0$ for all $u\in M$. Indeed, suppose on the contrary that $ G'(u)=0$ for some $u\in M $. Let
$$ \alpha:=a\int_{\mathbb R^{3}}|\nabla u|^{2},\ \ \ \beta:=\int_{\mathbb R^{3}}|u|^{2},\ \ \ \gamma:=b\left(\int_{\mathbb R^{3}}|\nabla u|^{2}\right)^{2},\ \ \ \theta:=\int_{\mathbb R^{3}}|u|^{p+1},\ \ \ \xi:=\int_{\mathbb R^{3}}|u|^{6},$$
and $ I(u)\triangleq l>0.$ The equation $ G'(u)=0$ can be written as
$$ -3\left(a+2b\int_{\mathbb R^{3}}|\nabla u|^{2}\right)\Delta u+5u-(p+4)|u|^{p-1}u-9|u|^{4}u=0, \eqno(2.7)$$
and $u$ satisfies the following Poho\u{z}aev identity
$$\frac{3}{2}a\int_{\mathbb R^{3}}|\nabla u|^{2}+\frac{15}{2}\int_{\mathbb R^{3}}|u|^{2}+3b\left(\int_{\mathbb R^{3}}|\nabla u|^{2}\right)^{2}-\frac{3(p+4)}{p+1}\int_{\mathbb R^{3}}|u|^{p+1}-\frac{9}{2}\int_{\mathbb R^{3}}|u|^{6}=0,$$
we then obtain
$$\left\{\aligned &\frac{1}{2}\alpha+\frac{1}{2}\beta+\frac{1}{4}\gamma-\frac{1}{p+1}\theta-\frac{1}{6}\xi=l,\\
&\frac{3}{2}\alpha+\frac{5}{2}\beta+\frac{3}{2}\gamma-\frac{p+4}{p+1}\theta-\frac{3}{2}\xi=0,\\
&3\alpha+5\beta+6\gamma-(p+4)\theta-9\xi=0,\\
&\frac{3}{2}\alpha+\frac{15}{2}\beta+3\gamma-\frac{3(p+4)}{p+1}\theta-\frac{9}{2}\xi=0. \endaligned \right. $$
The first equation comes from the fact that $ u$ minimizes $ I(u)=l$. The second one holds since $ G(u)=0.$ The third one follows from $ G'(u)u=0.$ The fourth one is the Poho\u{z}aev equality applied to (2.7). It can be checked out that the above system admits one unique solution for $ (\alpha,\beta,\gamma,\theta)$ depending on $ \xi, l, p.$ In particular, $ \theta$ is given by
$$ \theta=-\frac{10(p+1)}{3p(p-1)}(\xi+6l)<0,$$
which is contradicted with the definition of $\theta.$ So $ G'(u)\neq0$ for any $u\in M$ and by the Implicit Function theorem, $M$ is a $C^{1}$ manifold.

{\bf Step 4.} Every critical point of $ I|_{M}$ is a critical point of $I$ in $ H^{1}(\mathbb R^{3}).$

If $u$ is a critical point of $I|_{M}$, i.e., $u\in M$ and $(I|_{M})'(u)=0.$ Thanks to the Lagrange multiplier rule, there exists $ \rho\in \mathbb R$ such that $ I'(u)=\rho G'(u).$ We prove that $ \rho=0.$ As above, the equation $ I'(u)=\rho G'(u)$ can be written as
$$ \aligned &-\left(a+b\int_{\mathbb R^{3}}|\nabla u|^{2}\right)\Delta u+u-|u|^{p-1}u-|u|^{4}u\\
=&\ \rho\left(-3\left(a+2b\int_{\mathbb R^{3}}|\nabla u|^{2}\right)\Delta u+5u-(p+4)|u|^{p-1}u-9|u|^{4}u\right). \endaligned $$
So $ u$ can solve the equation
$$ -(3\rho-1)a\Delta u-(6\rho-1)b\int_{\mathbb R^{3}}|\nabla u|^{2}\Delta u+(5\rho-1)u-((p+4)\rho-1)|u|^{p-1}u-(9\rho-1)|u|^{4}u=0,$$
and $u$ satisfies the following Poho\u{z}aev identity
$$\aligned&\frac{3\rho-1}{2}a\int_{\mathbb R^{3}}|\nabla u|^{2}+\frac{3(5\rho-1)}{2}\int_{\mathbb R^{3}}|u|^{2}+\frac{6\rho-1}{2}b\left(\int_{\mathbb R^{3}}|\nabla u|^{2}\right)^{2}-\frac{3((p+4)\rho-1)}{p+1}\int_{\mathbb R^{3}}|u|^{p+1}\\
&-\frac{9\rho-1}{2}\int_{\mathbb R^{3}}|u|^{6}=0.\endaligned$$
By the definition of $\alpha,\beta,\gamma,\theta,\xi $ in {\bf Step 3}, we have
$$\left\{\aligned &\frac{1}{2}\alpha+\frac{1}{2}\beta+\frac{1}{4}\gamma-\frac{1}{p+1}\theta=l+\frac{1}{6}\xi,\\
&\frac{3}{2}\alpha+\frac{5}{2}\beta+\frac{3}{2}\gamma-\frac{p+4}{p+1}\theta=\frac{3}{2}\xi,\\
&(3\rho-1)\alpha+(5\rho-1)\beta+(6\rho-1)\gamma-((p+4)\rho-1)\theta=(9\rho-1)\xi,\\
&\frac{3\rho-1}{2}\alpha+\frac{3(5\rho-1)}{2}\beta+\frac{6\rho-1}{2}\gamma-\frac{3((p+4)\rho-1)}{p+1}\theta=\frac{9\rho-1}{2}\xi. \endaligned \right.\eqno(2.8) $$
The coefficient matrix of (2.8) is
$$B=\left( \begin{matrix}
\frac{1}{2}&\frac{1}{2}&\frac{1}{4}&-\frac{1}{p+1}\\
\frac{3}{2}& \frac{5}{2} & \frac{3}{2}&-\frac{p+4}{p+1}\\
3\rho-1 &5\rho-1&6\rho-1& -((p+4)\rho-1)\\
\frac{3\rho-1}{2}&\frac{3(5\rho-1)}{2}&\frac{6\rho-1}{2}&-\frac{3((p+4)\rho-1)}{p+1}\end{matrix}
\right) $$
and its determinant is
$$ \hbox{det}\ B=\frac{\rho(p-1)(2p-1-9p\rho)}{8(p+1)}.$$
We deduce immediately that
$$ \hbox{det}\ B=0\Longleftrightarrow\rho=0,\ \rho=\frac{2p-1}{9p},\ p=1.$$
Notice that $\rho$ must be equal to zero by excluding the other two possibilities:\

$(1)$ If $\rho\neq0,\rho\neq\frac{2p-1}{9p},$ the linear system (2.8) has a unique solution. Thus we deduce that
$$ \beta=-\frac{18l(p-5)((p+4)\rho-1)}{(p-1)(2p-1-9p\rho)},\ \ \ \ \ \ \theta=\frac{36l(p+1)(5\rho-1)}{(p-1)(2p-1-9p\rho)}.$$
Since $ p\in(2,5), \ \frac{1}{6}<\frac{2p-1}{9p}<\frac{1}{5},$ which shows that $\theta\leq0$ for $ \rho\in[\frac{1}{5},+\infty)\cup(-\infty,\frac{2p-1}{9p})$ and $ \beta<0$ for $ \rho\in(\frac{2p-1}{9p},\frac{1}{5}),$ however, this is impossible since $ \theta,\ \beta>0$.

$ (2)$ If $ \rho=\frac{2p-1}{9p},$ the latter two equations in (2.8) are as follows
$$\left\{\aligned &-\frac{p+1}{3p}\alpha+\frac{p-5}{9p}\beta+\frac{p-2}{3p}\gamma-\frac{2(p+1)(p-2)}{9p}\theta=0,\\
&-\frac{p+1}{6p}\alpha+\frac{p-5}{6p}\beta+\frac{p-2}{6p}\gamma-\frac{2(p-2)}{3p}\theta=0, \endaligned \right. $$
then $$ \beta+(p-2)\theta=0,$$
which is also impossible since $ p>2$, $ \beta,\ \theta>0$. Therefore $\rho=0$ and $ I'(u)=0$ in $ H^{1}(\mathbb R^{3}).$
$\hfill\square$

\vskip4pt

{\it\noindent {\bf Lemma 2.6.}  $c\triangleq c_{1}=c_{2}=c_{3}>0,$ where
$$c_{1}:=\inf_{\gamma\in \Gamma}\max_{t\in[0,1]}I(\gamma(t)),\ \ \ \ \ \ \ \ \ \ \ \ \ \ \ \ \ \ \ \ \ $$
$$c_{2}:=\inf_{u\in  H^{1}(\mathbb R^{3})\backslash\{0\}}\max_{t>0}I(u_{t}),\ \ \ \ \ \ \ \ \ \ \ \ \ \ $$
$$c_{3}:=\inf_{u\in M}I(u),\ \ \ \ \ \ \ \ \  \ \ \ \ \ \ \ \ \ \ \ \ \ \ \ \ \  \ \ \ \ $$
 $u_{t}$ is given in Lemma 2.2 and
$$\Gamma=\{\gamma\in C([0,1],H^{1}(\mathbb R^{3}))\ |\ \gamma(0)=0,\ I(\gamma(1))\leq0,\ \gamma(1)\neq0\}.$$}
\noindent{\it Proof.} The proof is similar to the proof of ($\cite{31}$, Lemma 2.8), where $ M$ was the manifold of Nehari-Poho\u{z}aev type. The proof consists three steps.

{\bf Step 1.} $c_{3}>0.$ It could be directly deduced by {\bf Step 2} in the proof of Lemma 2.5.

{\bf Step 2.} $c_{2}=c_{3}.$ Using Lemma 2.4, for every $u\in H^{1}(\mathbb R^{3})\backslash\{0\},$ there exists a unique $u_{\hat{t}}\in M$ such that $ I(u_{\hat{t}})=\max\limits_{t>0}I(u_{t}).$

{\bf Step 3.} $c_{2}\geq c_{1}\geq c_{3}.$ For every $ \gamma\in \Gamma,$ we show that $ \gamma([0,1]\cap M)\neq\emptyset.$ In fact, combing  {\bf Step 1} with {\bf Step 2} in the proof of Lemma 2.5, we see that $\gamma$ crosses $ M$ since $ \gamma(0)=0,\ I(\gamma(1))\leq0$ and $ \gamma(1)\neq0.$ Therefore $ \max\limits_{t\in[0,1]}I(\gamma(t))\geq\inf\limits_{u\in M}I(u)=c_{3},$
which implies that $ c_{1}\geq c_{3}.$
Furthermore, for $u\in (H^{1}(\mathbb R^{3})\cap C(\mathbb R^{3}))\backslash\{0\},$ by Lemma 2.2, $ I(u_{t_{0}})<0$ for $ t_{0}$ large enough.
Let $$ \bar{\gamma}(t)\triangleq\left\{\aligned &u_{tt_{0}},&t>0,\\
&0,&t=0, \endaligned \right. $$
which implies that $ \bar{\gamma}\in \Gamma.$ Then
$$ \max_{t>0}I(u_{tt_{0}})\geq\max_{t\in[0,1]}I(u_{tt_{0}})\geq\inf_{\gamma\in \Gamma}\max_{t\in[0,1]}I(\gamma(t))=c_{1}.$$
Thus we conclude $ c_{2}\geq c_{1}.$  $\hfill\square$

If $u\in M$ such that $I(u)=c$, then $u$ is a ground state solution of Nehari-Poho\u{z}aev type to (1.7) by Lemma 2.5 and Lemma 2.6. Therefore we look for critical points of $I|_{M}$.

The following concentration-compactness principle is due to P. L. Lions' results.

\noindent {\bf Lemma 2.7.} (\cite[Lemma 1.1]{33}) {\it Let $\{\rho_{n}\}$ be a sequence of nonnegative $L^{1}$ functions on $ \mathbb R^{N}$ satisfying $ \int_{\mathbb R^{N}}\rho_{n}=s, $ where $ s>0$ is fixed. There exists a subsequence, still denoted by $\{\rho_{n}\}$ satisfying one of the following three possibilities:\\
(i)(Vanishing) for all $ R>0,$ it holds
$$ \lim_{n\rightarrow\infty}\sup_{y\in \mathbb R^{N}}\int_{B_{R}(y)}\rho_{n}=0;$$
(ii)(Compactness) there exists $ \{y_{n}\}\subset \mathbb R^{N}$ such that, for any $ \varepsilon>0,$ there exists $ R>0$ satisfying
$$ \liminf_{n\rightarrow\infty}\int_{B_{R}(y_{n})}\rho_{n}\geq s-\varepsilon;$$
(iii)(Dichotomy) there exist $\alpha\in(0,s) $ and $ \{y_{n}\}\subset \mathbb R^{N}$ such that for $ \varepsilon>0, \ \hbox{there exists} \ R>0,$ for all $ r\geq R$ and $ r'\geq R,$ it holds
$$ \limsup_{n\rightarrow\infty}\left|\alpha-\int_{B_{R}(y_{n})}\rho_{n}\right|+\left|(s-\alpha)-\int_{\mathbb R^{N}\backslash B_{r'}(y_{n})}\rho_{n}\right|<\varepsilon.$$ }

\noindent {\bf Lemma 2.8.} {\it Let $ r>0, \ s\in[2,2^{*}].$ If $ \{u_{n}\}$ is bounded in $ H^{1}(\mathbb R^{N})$ and
$$ \lim_{n\rightarrow\infty}\sup_{y\in \mathbb R^{N}}\int_{B_{r}(y)}|u_{n}|^{p}=0,$$
then we have $ u_{n}\rightarrow0$ in $ L^{q}(\mathbb R^{N})$ for $ q\in(2,2^{*}).$ Moreover, if $p = 2^{*}$, $u_{n}\rightarrow0$ in $L^{2^{*}}(\mathbb R^{N})$. Here $2^{*}= \frac{2N}{N-2}$ if $N\geq3$ and $2^{*}=\infty$ if $N =1,2.$}

For a proof of the first part of this lemma see \cite{34} and for the second part see \cite{0} (and \cite{51} for similar ideas).

\noindent {\bf Lemma 2.9.} (\cite[Lemma 1.32]{34}) {\it Let $\Omega$ be an open subset of $\mathbb R^{N}$ and let $\{u_{n}\}\subset L^{p}(\Omega), 1\leq p<\infty.$ If $\{u_{n}\}$ is bounded in $L^{p}(\Omega)$ and $u_{n}\rightarrow u$ a.e. on $\Omega$, then
$$\lim_{n\rightarrow\infty}(|u_{n}|_{p}^{p}-|u_{n}-u|_{p}^{p})=|u|_{p}^{p}.$$}

\section{ Existence of positive ground state of Nehari-Poho\u{z}aev type}
 In this section, we prove the existence of positive ground state solution of Nehari-Poho\u{z}aev type to $ (K_{\varepsilon}).$

{\it\noindent {\bf Lemma 3.1.} Let $\{u_{n}\}\subset M $ be a minimizing sequence for $ c,$ which was given in Lemma 2.6, then there exists $ \{y_{n}\}\subset \mathbb R^{3}$ such that for $ \varepsilon>0,$ there exists $ R>0$ such that
$$ \int_{\mathbb R^{3}\backslash B_{R}(y_{n})}(a|\nabla u_{n}|^{2}+|u_{n}|^{2})\leq \varepsilon.$$}

\noindent {\it Proof.}  Let $ \{u_{n}\}\subset M$ be a minimizing sequence of $I$ at the level $c,$ satisfying
$$ \lim_{n\rightarrow+\infty}I(u_{n})=c>0. \eqno(3.1)$$
We define $ \Psi:H^{1}(\mathbb R^{3})\rightarrow \mathbb R$ as
$$\Psi(u)=\int_{\mathbb R^{3}}\left(\frac{a}{4}|\nabla u|^{2}+\frac{1}{12}|u|^{2}+\frac{p-2}{6(p+1)}|u|^{p+1}+\frac{1}{12}|u|^{6}\right). \eqno(3.2)$$
In fact, $\Psi(u)=I(u)-\frac{1}{6}G(u).$ For any $ u\in M,$ we have that $ I(u)=\Psi(u)\geq0$ and then $ \lim\limits_{n\rightarrow\infty}\Psi(u_{n})=c.$ By (3.2), we see that for $n$ sufficiently large,
$c+1+\|u_{n}\|_{H}\geq \Psi(u_{n})\geq \frac{1}{12}\|u_{n}\|_{H}^{2}, $
then $\{u_{n}\}$ is bounded in $ H^{1}(\mathbb R^{3}).$ Up to a subsequence, we can assume that there exists $u\in H^{1}(\mathbb R^{3})$ such that $u_{n}\rightharpoonup u\ \hbox{weakly in} \ H^{1}(\mathbb R^{3}).$ By the Rellich theorem, $u_{n}\rightarrow u\ \hbox{ in} \ L_{loc}^{q}(\mathbb R^{3})$, $q\in[1,6).$
Let
$$ \rho_{n}\triangleq \frac{a}{4}|\nabla u_{n}|^{2}+\frac{1}{12}|u_{n}|^{2}+\frac{p-2}{6(p+1)}|u_{n}|^{p+1}+\frac{1}{12}|u_{n}|^{6},\ \ x\in \mathbb R^{3},\eqno(3.3)$$
then $ \{\rho_{n}\}$ is a sequence of nonnegative $L^{1}$ function on $ \mathbb R^{3}$ satisfying
$$ \int_{\mathbb R^{3}}\rho_{n}=\Psi(u_{n})\rightarrow c>0.$$
By Lemma 2.7, there are three cases:

(1) Vanishing: for all $ R>0,$
$$ \lim_{n\rightarrow+\infty}\sup_{y\in \mathbb R^{3}}\int_{B_{R}(y_{n})}\rho_{n}=0.$$

(2) Compactness: there exists $ \{y_{n}\}\subset \mathbb R^{3}$ such that for $\varepsilon>0,$ there exists $ R>0$ satisfying
$$\liminf_{n\rightarrow+\infty}\int_{B_{R}(y_{n})}\rho_{n}\geq c-\varepsilon.$$

(3) Dichotomy: there exists $\alpha\in(0,c)$ and $\{y_{n}\}\subset \mathbb R^{3}$ such that for $ \varepsilon>0,$ there exists $R>0$ satisfying
$$ \limsup_{n\rightarrow\infty}\left|\alpha-\int_{B_{R}(y_{n})}\rho_{n}\right|+\left|(c-\alpha)-\int_{\mathbb R^{3}\backslash B_{2R}(y_{n})}\rho_{n}\right|<\varepsilon.$$
Now we show that compactness holds for the sequence $\{\rho_{n}\}$ defined in (3.3).

(i) Vanishing cannot occur.\\
\indent If $\{\rho_{n}\}$ vanishes, then $\{|u_{n}|^{p+1}+u_{n}^{6}\}$ are vanishes, i.e., there exists $R>0$ such that
$$\lim_{n\rightarrow+\infty}\sup_{y\in {\mathbb R^{3}}}\int_{B_{R}(y)}(|u_{n}|^{p+1}+|u_{n}|^{6})=0.$$
In view of Lemma 2.8, one has $u_{n}\rightarrow0 $ in $L^{q}({\mathbb R^{3}})$, $q\in(2,6]$. Since $\{u_{n}\}\subset M$ and $p\in(2,5)$,
$$\aligned 0<c&\leq I(u_{n})\\
&=\ I(u_{n})-G(u_{n})\\
&=-a\int_{\mathbb R^{3}}|\nabla u_{n}|^{2}-2\int_{\mathbb R^{3}}|u_{n}|^{2}-\frac{5}{4}b\left(\int_{\mathbb R^{3}}|\nabla u_{n}|^{2}\right)^{2}+\frac{p+3}{p+1}\int_{\mathbb R^{3}}|u_{n}|^{p+1}+\frac{4}{3}\int_{\mathbb R^{3}}|u_{n}|^{6}\\
&<\frac{p+3}{p+1}\int_{\mathbb R^{3}}|u_{n}|^{p+1}+\frac{4}{3}\int_{\mathbb R^{3}}|u_{n}|^{6}\\
&\rightarrow0,\endaligned$$
which is a contradiction.

(ii) Dichotomy cannot occur.

If there exist $ \alpha\in(0,c)$ and $\{y_{n}\}\subset \mathbb R^{3} $ such that for $ \varepsilon_{n}\rightarrow0,$ we can choose $\{R_{n}\}\subset \mathbb R^{+}$ with $R_{n}\rightarrow\infty$ satisfying
$$ \limsup_{n\rightarrow+\infty}\left|\alpha-\int_{B_{R_{n}}(y_{n})}\rho_{n}\right|+\left|(c_{\varepsilon}-\alpha)-\int_{\mathbb R^{3}\backslash B_{2R_{n}}(y_{n})}\rho_{n}\right|<\varepsilon_{n}. \eqno(3.4)$$
Let $ \xi:\mathbb R^{+}\rightarrow \mathbb R^{+}$ be a cut-off function satisfying $ 0\leq\xi\leq1, \ \xi(s)\equiv1 $ for $ s\leq1, \xi\equiv0$ for $ s\geq2$ and $ |\xi'(s)|\leq2,$
$$ v_{n}(x)\triangleq\xi\left( \frac{|x-y_{n}|}{R_{n}}\right)u_{n}(x), \ \ \ w_{n}(x)\triangleq\left(1-\xi\left( \frac{|x-y_{n}|}{R_{n}}\right)\right)u_{n}(x).$$
From (3.4) we deduce immediately that $$\liminf_{n\rightarrow\infty}\Psi(v_{n})\geq\alpha, \ \ \ \ \ \ \liminf_{n\rightarrow\infty}\Psi(w_{n})\geq c-\alpha. \eqno(3.5)$$
Denote $ \Omega_{n}:=B_{2R_{n}(y_{n})}\backslash B_{R_{n}(y_{n})},$ passing to the limit as $ n\rightarrow+\infty$,
$$ \int_{\Omega_{n}}\left(\frac{a}{4}|\nabla u_{n}|^{2}+\frac{1}{12}|u_{n}|^{2}+\frac{p-2}{6(p+1)}|u_{n}|^{p+1}+\frac{1}{12}|u_{n}|^{6}\right)=\int_{\Omega_{n}}\rho_{n}\rightarrow0.$$
Then,
$$ \int_{\Omega_{n}}\left(a|\nabla u_{n}|^{2}+|u_{n}|^{2}\right)\rightarrow0, \ \ \int_{\Omega_{n}}|u_{n}|^{p+1}\rightarrow0,\ \  \ \int_{\Omega_{n}}|u_{n}|^{6}\rightarrow0\ \ \hbox{as}\ \ n\rightarrow+\infty.$$
and
$$ \int_{\Omega_{n}}\left(a|\nabla v_{n}|^{2}+|v_{n}|^{2}\right)\rightarrow0, \ \ \int_{\Omega_{n}}\left(a|\nabla w_{n}|^{2}+|w_{n}|^{2}\right)\rightarrow0\ \ \hbox{as}\ \ n\rightarrow+\infty.$$
one has
$$\aligned & a\int_{\mathbb R^{3}}|\nabla u_{n}|^{2}=a\int_{\mathbb R^{3}}|\nabla v_{n}|^{2}+a\int_{\mathbb R^{3}}|\nabla w_{n}|^{2}+o(1),\\
&\int_{\mathbb R^{3}}|u_{n}|^{2}=\int_{\mathbb R^{3}}|v_{n}|^{2}+\int_{\mathbb R^{3}}|w_{n}|^{2}+o(1),\\
&\int_{\mathbb R^{3}}|u_{n}|^{p+1}=\int_{\mathbb R^{3}}|v_{n}|^{p+1}+\int_{\mathbb R^{3}}|w_{n}|^{p+1}+o(1),\\
&\int_{\mathbb R^{3}}|u_{n}|^{6}=\int_{\mathbb R^{3}}|v_{n}|^{6}+\int_{\mathbb R^{3}}|w_{n}|^{6}+o(1).\endaligned\eqno(3.6) $$
Moreover,
$$ \aligned \left(\int_{\mathbb R^{3}}|\nabla u_{n}|^{2}\right)^{2}&=\left(\int_{\mathbb R^{3}}|\nabla v_{n}|^{2}+\int_{\mathbb R^{3}}|\nabla w_{n}|^{2}+o(1)\right)^{2}\\
&\geq \left(\int_{\mathbb R^{3}}|\nabla v_{n}|^{2}\right)^{2}+\left(\int_{\mathbb R^{3}}|\nabla w_{n}|^{2}\right)^{2}+o(1). \endaligned \eqno(3.7)$$
By (3.5)(3.6),
$$ \Psi(u_{n})=\Psi(v_{n})+\Psi(w_{n})+o(1),$$
$$c=\lim_{n\rightarrow+\infty}\Psi(u_{n})\geq\liminf_{n\rightarrow+\infty}\Psi(v_{n})+\liminf_{n\rightarrow+\infty}\Psi(w_{n})\geq\alpha+c-\alpha=c,$$
therefore $$ \lim_{n\rightarrow+\infty}\Psi(v_{n})=\alpha,\ \ \ \lim_{n\rightarrow+\infty}\Psi(w_{n})=c-\alpha. \eqno(3.8)$$
By (3.5), (3.6) and $u_{n}\in M$, we obtain
$$ 0=G(u_{n})\geq G(v_{n})+G(w_{n})+o(1). \eqno(3.9)$$

Next, we distinguish two cases:

{\bf Case 1.} Up to a subsequence, we may assume that either $ G(v_{n})\leq0$ or $ G(w_{n})\leq0.$

Without loss of generality, we assume $ G(v_{n})\leq0,$ then
$$\frac{3}{2}a\int_{\mathbb R^{3}}|\nabla v_{n}|^{2}+\frac{5}{2}\int_{\mathbb R^{3}}|v_{n}|^{2}+\frac{3}{2}b\left(\int_{\mathbb R^{3}}|\nabla v_{n}|^{2}\right)^{2}-\frac{p+4}{p+1}\int_{\mathbb R^{3}}|v_{n}|^{p+1}-\frac{3}{2}\int_{\mathbb R^{3}}|v_{n}|^{6}\leq0. \eqno(3.10)$$
For every $n$, there exists $ t_{n}>0$ such that $ (v_{n})_{t_{n}}\in M$ by Lemma 2.4 and then $ G((v_{n})_{t_{n}})=0, $ i.e.,
$$\aligned  &\frac{3}{2}at_{n}^{3}\int_{\mathbb R^{3}}|\nabla v_{n}|^{2}+\frac{5}{2}t_{n}^{5}\int_{\mathbb R^{3}}|v_{n}|^{2}+\frac{3}{2}bt_{n}^{6}\left(\int_{\mathbb R^{3}}|\nabla v_{n}|^{2}\right)^{2}-\frac{p+4}{p+1}t_{n}^{p+4}\int_{\mathbb R^{3}}|v_{n}|^{p+1}\\
&-\frac{3}{2}t_{n}^{9}\int_{\mathbb R^{3}}|v_{n}|^{6}=0. \endaligned \eqno(3.11)$$
Notice that (3.10) and (3.11) imply that
$$ \aligned &\frac{3}{2}a(t_{n}^{p+1}-1)\int_{\mathbb R^{3}}|\nabla v_{n}|^{2}+\frac{5}{2}(t_{n}^{p+1}-t_{n}^{2})\int_{\mathbb R^{3}}|v_{n}|^{2}+\frac{3}{2}b(t_{n}^{p+1}-t_{n}^{3})\left(\int_{\mathbb R^{3}}|\nabla v_{n}|^{2}\right)^{2}\\
&-\frac{3}{2}(t_{n}^{p+1}-t_{n}^{6})\int_{\mathbb R^{3}}|v_{n}|^{6}\leq0, \endaligned$$
so that $ t_{n}\leq1.$ We would have
$$ c\leq I((v_{n})_{t_{n}})=\Psi((v_{n})_{t_{n}})\leq\Psi(v_{n})\rightarrow\alpha<c, \eqno(3.12)$$
a contradiction.

{\bf Case 2.} Up to a subsequence, we may suppose that $ G(v_{n})>0$ and $G(w_{n})>0.$

We see that $ G(v_{n})\rightarrow0$ and $ G(w_{n})\rightarrow0$ as $ n\rightarrow+\infty$ by (3.9). For $ t_{n}$ given in Case 1, if $\limsup\limits_{n\rightarrow+\infty}t_{n}\leq1,$ then we can get the same contradiction as (3.12). Assume now that $ \lim\limits_{n\rightarrow+\infty}t_{n}=t_{0}>1,$ by (3.11), we have
$$\aligned  G(v_{n})&=\frac{3}{2}a\int_{\mathbb R^{3}}|\nabla v_{n}|^{2}+\frac{5}{2}\int_{\mathbb R^{3}}|v_{n}|^{2}+\frac{3}{2}b\left(\int_{\mathbb R^{3}}|\nabla v_{n}|^{2}\right)^{2}-\frac{p+4}{p+1}\int_{\mathbb R^{3}}|v_{n}|^{p+1}-\frac{3}{2}\int_{\mathbb R^{3}}|v_{n}|^{6}\\
&=\frac{3}{2}a\left(1-\frac{1}{t_{n}^{p+1}}\right)\int_{\mathbb R^{3}}|\nabla v_{n}|^{2}+\frac{5}{2}\left(1-\frac{1}{t_{n}^{p-1}}\right)\int_{\mathbb R^{3}}|v_{n}|^{2}+\frac{3}{2}b\left(1-\frac{1}{t_{n}^{p-2}}\right)\left(\int_{\mathbb R^{3}}|\nabla v_{n}|^{2}\right)^{2}\\
&\ \ \ -\frac{3}{2}\left(1-\frac{1}{t_{n}^{p-5}}\right)\int_{\mathbb R^{3}}|v_{n}|^{6}. \endaligned $$
Then $ v_{n}\rightarrow0$ in $ H^{1}(\mathbb R^{3})$ since $ G(v_{n})\rightarrow0$ as $n\rightarrow\infty$, which contradicts to (3.8) since $ \alpha>0,$ which shows that dichotomy cannot occur.\\
\indent Therefore, the sequence $\{\rho_{n}\}$ is compact, i.e., there exists $\{y_{n}\}\subset \mathbb R^{3}$ such that for $ \varepsilon>0,$ there exists $ R>0$ satisfying
$$ \liminf_{n\rightarrow\infty}\int_{B_{R}(y_{n})}\left(\frac{a}{4}|\nabla u_{n}|^{2}+\frac{1}{12}|u_{n}|^{2}+\frac{p-2}{6(p+1)}|u_{n}|^{p+1}+\frac{1}{12}|u_{n}|^{6}\right)\geq c-\varepsilon.$$
It follows from $\lim\limits_{n\rightarrow+\infty}\Psi(u_{n})=c$ that $ \int_{\mathbb R^{3}\backslash B_{R}(y_{n})}a|\nabla u_{n}|^{2}+|u_{n}|^{2}\leq\varepsilon.$ $\hfill\square$

\vskip4pt

{\it\noindent {\bf Proof of Theorem 1.6} } Let $ \{u_{n}\}\subset M$ be a minimizing sequence for $ c$ (see Lemma 2.6), then by Lemma 3.1, there exists $\{y_{n}\}\subset \mathbb R^{3}$ such that for $ \varepsilon>0,$ there exists $ R>0$ satisfying
$$ \int_{\mathbb R^{3}\backslash B_{R}(y_{n})}\left(a|\nabla u_{n}|^{2}+|u_{n}|^{2}\right)\leq \varepsilon. \eqno(3.13)$$
Let $ \bar{u}_{n}(\cdot)={u}_{n}(\cdot-y_{n})\in H^{1}(\mathbb R^{3}),$ then $ \bar{u}_{n}\in M$,
$$ \int_{\mathbb R^{3}\backslash B_{R}(0)}\left(a|\nabla \bar{u}_{n}|^{2}+|\bar{u}_{n}|^{2}\right)\leq \varepsilon. \eqno(3.14)$$
Since $\{\bar{u}_{n}\}$ is bounded in $ H^{1}(\mathbb R^{3}),$ up to a subsequence, we may assume that there exists $ \bar{u}\in H^{1}(\mathbb R^{3})$ such that
$$ \left\{\aligned &\bar{u}_{n}\rightharpoonup \bar{u},\ \ \ \ \ \ \ \ \ \ \ \hbox{in} \ \ H^{1}(\mathbb R^{3}),\\
&\bar{u}_{n}\rightarrow \bar{u},\ \ \ \ \ \ \ \ \ \ \ \hbox{in} \ \ L_{loc}^{s}(\mathbb R^{3}), \ \ s\in[2,6),\\
&\bar{u}_{n}(x)\rightarrow \bar{u}(x),\ \  \ \hbox{a.e.  in}\ \ \mathbb R^{3}, \endaligned  \right. \eqno(3.15)  $$
and
$$ \int_{\mathbb R^{3}\setminus B_{R}(0)}(a|\nabla \bar{u}|^{2}+|\bar{u}|^{2})\leq \varepsilon. \eqno(3.16)$$
By (3.14)-(3.16),  for $ s\in [2,6)$ and $ \varepsilon>0,$ there exists $ C>0$ such that
$$ \aligned \int_{\mathbb R^{3}}|\bar{u}_{n}-\bar{u}|^{s}&=\int_{B_{R}(0)}|\bar{u}_{n}-\bar{u}|^{s}+\int_{\mathbb R^{3}\setminus B_{R}(0)}|\bar{u}_{n}-\bar{u}|^{s}\\
&\leq \varepsilon+C(\|\bar{u}_{n}\|_{H^{1}\left(\mathbb R^{3}\setminus B_{R}(0)\right)}+\|\bar{u}\|_{H^{1}(\mathbb R^{3}\setminus B_{R}(0))})\\
&\leq(1+2C)\varepsilon.  \endaligned \eqno(3.17)$$
Therefore
$$ \bar{u}_{n}\rightarrow \bar{u} \ \  \hbox{ in}  \ L^{s}(\mathbb R^{3}), \ \ \ s\in [2,6).  \eqno(3.18)$$
Since $\bar{u}_{n}\in M,$ by Sobolev embedding theorem, we have
$$3a\int_{\mathbb R^{3}}|\nabla \bar{u}_{n}|^{2}+5\int_{\mathbb R^{3}}|\bar{u}_{n}|^{2}+3b\left(\int_{\mathbb R^{3}}|\nabla \bar{u}_{n}|^{2}\right)^{2}-\frac{2(p+4)} {p+1}\int_{\mathbb R^{3}}|\bar{u}_{n}|^{p+1}-3\int_{\mathbb R^{3}}|\bar{u}_{n}|^{6}=0,$$
thus,
$$\aligned 3\|\bar{u}_{n}\|_{H}^{2}&\leq\frac{2(p+4)} {p+1}\int_{\mathbb R^{3}}|\bar{u}_{n}|^{p+1}+3\int_{\mathbb R^{3}}|\bar{u}_{n}|^{6}\\
&=\frac{2(p+4)} {p+1}|\bar{u}_{n}|_{p+1}^{p+1}+3|\bar{u}_{n}|_{6}^{6}\\
&\leq C\|\bar{u}_{n}\|_{H}^{6}.\endaligned$$
So for some $C> 0$, we get that $\|\bar{u}_{n}\|_{H}\geq C>0,$ $|\bar{u}_{n}|_{p+1}\geq C>0$, which implies that $\bar{u}\neq0.$

We claim that $\bar{u}\in M$ and $ \bar{u}_{n}\rightarrow \bar{u}$ in $ H^{1}(\mathbb R^{3}).$ In fact, if $G(\bar{u})<0$, by Lemma 2.4, there exists $ 0<t_{0}<1$ such that $ \bar{u}_{t_{0}}\in M.$ Then $ G(\bar{u}_{t_{0}})=0$ and
$$ I(\bar{u}_{t_{0}})=\Psi(\bar{u}_{t_{0}})<\Psi(\bar{u})\leq\lim_{n\rightarrow+\infty}\Psi(\bar{u}_{n})=\lim_{n\rightarrow+\infty}I(\bar{u}_{n})=c,$$
which contradicts the definition of $c$, where $\Psi$ is given in (3.2). If $G(\bar{u})>0$, set $\xi_{n}:=\bar{u}_{n}-\bar{u}$, by Lemma 2.9, $(V_{2})$, we may obtain
\begin{equation}\label{eq3.7}
G(\bar{u}_{n})=G(\bar{u})+G(\xi_{n})+o_{n}(1).
\end{equation}
Then $\limsup\limits_{n\rightarrow\infty}G(\xi_{n})<0$.  By Lemma 2.4, there exists $t_{n}\in(0,1)$ such that $(\xi_{n})_{t_{n}}\in M.$ Furthermore, one has that $\limsup\limits_{n\rightarrow\infty}t_{n}<1$, otherwise, along a subsequence, $t_{n}\rightarrow1$ and hence $G(\xi_{n})=G((\xi_{n})_{t_{n}})+o_{n}(1) = o_{n}(1)$, a contradiction. It follows from $\bar{u}_{n}\in M$, (3.1), $(V_{2}) $  that
$$\aligned &\ c+o_{n}(1)\\
=&\ I( \bar{u}_{n} )-\frac{1}{6}G(\bar{u}_{n})\\
=&\ \int_{\mathbb R^{3}}\left(\frac{a}{4}|\nabla \bar{u}_{n}|^{2}+\frac{1}{12}|\bar{u}_{n}|^{2}+\frac{p-2}{6(p+1)}|\bar{u}_{n}|^{p+1}+\frac{1}{12}|\bar{u}_{n}|^{6}\right)\\
=&\ \int_{\mathbb R^{3}}\bigg(\frac{a}{4}(|\nabla \bar{u}_{n}|^{2}+|\nabla \xi_{n}|^{2})+\frac{1}{12}(|\bar{u}_{n}|^{2}+|\xi_{n}|^{2})+\frac{p-2}{6(p+1)}(|\bar{u}_{n}|^{p+1}+|\xi_{n}|^{p+1})\\
&\ +\frac{1}{12}(|\bar{u}_{n}|^{6}+|\xi_{n}|^{6})\bigg)\\
=&\ I( \bar{u})-\frac{1}{6}G(\bar{u})+I(\xi_{n})-\frac{1}{6}G(\xi_{n})\\
>&\ I( \bar{u} )-\frac{1}{6}G(\bar{u})+ I((\xi_{n}) _{t_{n}})-\frac{1}{6}G((\xi_{n})_{t_{n}})\\
=&\ I((\xi_{n})_{t})+\int_{\mathbb R^{3}}\left(\frac{a}{4}|\nabla \bar{u}|^{2}+\frac{1}{12}|\bar{u}|^{2}+\frac{p-2}{6(p+1)}|\bar{u}|^{p+1}+\frac{1}{12}|\bar{u}|^{6}\right)\\
\geq&\ c, \endaligned$$
which is also a contradiction. So $ \bar{u}\in M$. By Fatou's Lemma, we may obtain
$$\aligned c&\leq I(\bar{u})\\
&=I(\bar{u})-\frac{1}{6}G(\bar{u})\\
&=\int_{\mathbb R^{3}}\left(\frac{a}{4}|\nabla \bar{u}|^{2}+\frac{1}{12}|\bar{u}|^{2}+\frac{p-2}{6(p+1)}|\bar{u}|^{p+1}+\frac{1}{12}|\bar{u}|^{6}\right)\\
&\leq \liminf_{n\rightarrow\infty}\bigg(\frac{a}{4}|\nabla \bar{u}_{n}|^{2}+\frac{1}{12}|\bar{u}_{n}|^{2}+\frac{p-2}{6(p+1)}|\bar{u}_{n}|^{p+1}+\frac{1}{12}|\bar{u}_{n}|^{6}\bigg)\\
&=\liminf_{n\rightarrow\infty}\bigg(I(\bar{u}_{n})-\frac{1}{6}G(\bar{u}_{n})\bigg)\\
&=c.\endaligned$$
Therefore, $\bar{u}_{n}\rightarrow\bar{u}$ in $ H^{1}(\mathbb R^{3}).$

We deduce that $ \bar{u}\in M$ and $ I(\bar{u})=c,$ i.e., $ I|_{M}$ attains its minimum at $ \bar{u},$ then $\bar{u}$ is a nontrivial critical point of $I|_{M},$ hence by Lemma 2.5, we see that $ \bar{u}$ is a ground state solution of Nehari-Poho\u{z}aev type of (1.6). Note that $ |\bar{u}| $ is also a ground state solution of Nehari-Poho\u{z}aev type of (1.6) since the functional $I$ and the manifold $ M$ are symmetric, so that we can assume that $\bar{u}\geq0.$ The strong maximum principle and standard arguments  $\cite{36, 37}$ imply that $\bar{u}(x)>0 $ for all $ x\in \mathbb R^{3}$ and the Theorem 1.6 is proved.    $\hfill\square$

If $(V_{1})$-$(V_{4})$, $(F_{1})$-$(F_{5})$ hold, we apply Lemma 1.5 to prove Theorem 1.1 (1). Set $J=[\delta,1]$, where $\delta\in(0,1)$. We investigate a family of functionals on $E_{\varepsilon}$:
$$ I_{\varepsilon,\sigma}(u)=\frac{1}{2}\|u\|_{\varepsilon}^{2}+\frac{b}{4}\left(\int_{\mathbb R^{3}}|\nabla u|^{2}\right)^{2}-\frac{\sigma}{p+1}\int_{\mathbb R^{3}}A(\varepsilon x)|u|^{p+1}-\frac{\sigma}{6}\int_{\mathbb R^{3}}B(\varepsilon x)|u|^{6}.\eqno(3.20)$$
Then let $I_{\varepsilon,\sigma}(u)=C(u)-\sigma D(u),$ where
$$ C(u):=\frac{1}{2}\|u\|_{\varepsilon}^{2}+\frac{b}{4}\left(\int_{\mathbb R^{3}}|\nabla u|^{2}\right)^{2}\rightarrow+\infty\ \ \hbox{as} \ \ \|u\|_{\varepsilon}\rightarrow+\infty,\ $$
$$ D(u):=\frac{1}{p+1}\int_{\mathbb R^{3}}A(\varepsilon x)|u|^{p+1}+\frac{1}{6}\int_{\mathbb R^{3}}B(\varepsilon x)|u|^{6}\geq0.\ \ \ \ \ \ \ \ \ $$

{\it\noindent {\bf Lemma 3.4.} If $(V_{2})$, $(F_{3})$-$(F_{5})$ hold and $ p\in(2,5).$ Then\\
$(1)$ there exists $ v\in E_{\varepsilon}\backslash\{0\}$ such that $ I_{\varepsilon,\sigma}(v)\leq0,$ for $ \sigma\in[\delta,1];$ \\
$(2)$ $ c_{\sigma}:=\inf\limits_{\gamma \in \Gamma}\max\limits_{s\in[0,1]}I_{\varepsilon,\sigma}(\gamma(s))>\max\{I_{\varepsilon,\sigma}(0), \ I_{\varepsilon,\sigma}(v)\}$ for $\sigma\in[\delta,1],$ where\\ \indent $\Gamma=\{\gamma\in C([0,1], E_{\varepsilon})|\ \gamma(0)=0,\ \gamma(1)=v\}.$
}

\noindent {\it Proof.} The proof is standard and we omit it here. $ \hfill\square$

Lemma 3.4 means that $I_{\varepsilon,\sigma}(u)$ satisfies the assumptions of Lemma 1.5 with $X=E_{\varepsilon}$ and $ \Phi_{\sigma}=I_{\varepsilon,\sigma}.$ We then obtain immediately, for a.e. $ \sigma\in [\delta,1],$ there exists a bounded sequence $ \{u_{n}\}\subset E_{\varepsilon}$ such that
$ I_{\varepsilon,\sigma}(u_{n})\rightarrow c_{\sigma},\ I'_{\varepsilon,\sigma}(u_{n})\rightarrow0 \ \ \hbox{in}\ \ E_{\varepsilon}.$

\noindent {\bf Lemma 3.5.} (\cite{27}, Lemma 2.3) {\it Under the assumptions of Lemma 1.5, the map $ \sigma\rightarrow c_{\sigma}$ is non-increasing and left continuous. }

By Theorem 1.6, we conclude that for any $\sigma\in[\delta,1]$, the associated limit problem (1.5) has a positive ground state solution of Nehari-Poho\u{z}aev type in $E_{\varepsilon}$, i.e., for any $\sigma\in[\delta,1]$,
$$ m_{\sigma}^{\infty}:=\inf_{u\in M_{\sigma}^{\infty}}I_{\sigma}^{\infty}(u) \eqno(3.21)$$
is achieved at some $u_{\sigma}^{\infty}\in M_{\sigma}^{\infty}:=\{u\in E_{\varepsilon}\backslash\{0\}\ |\ G_{\sigma}^{\infty}(u)=0\},$ and $(I_{\sigma}^{\infty})'(u_{\sigma}^{\infty})=0,$
where
$$G_{\sigma}^{\infty}:=\frac{3a}{2}\int_{\mathbb R^{3}}|\nabla u|^{2}+\frac{5}{2}\int_{\mathbb R^{3}}V_{\infty}|u|^{2}+\frac{3b}{2}\left(\int_{\mathbb R^{3}}|\nabla u|^{2}\right)^{2}-\frac{(p+4)\sigma}{p+1}\int_{\mathbb R^{3}}A_{\infty}|u|^{p+1}-\frac{3\sigma}{2}\int_{\mathbb R^{3}}B_{\infty}|u|^{6}. \eqno(3.22)$$

{\it\noindent {\bf Lemma 3.6.} Suppose that $ (V_{1})$, $(V_{2})$, $(F_{3})$, $(F_{4})$ hold, $2<p<5.$ Then $ c_{\sigma}<m_{\sigma}^{\infty}$ for $ \sigma\in[\delta,1].$ }\\
\noindent {\it Proof.} Let $ u_{\sigma}^{\infty}$ be the minimizer of $m_{\sigma}^{\infty}$, by Lemma 2.4, $I_{\sigma}^{\infty}(u_{\sigma}^{\infty})=\max_{t>0}I_{\sigma}^{\infty}(tu_{\sigma}^{\infty}(t^{-1}x)).$
Then,  for $\sigma\in[\delta,1],$
$$ c_{\sigma}\leq \max_{t>0}I_{\varepsilon,\sigma}(tu_{\sigma}^{\infty}(t^{-1}x))<\max_{t>0}I_{\sigma}^{\infty}(tu_{\sigma}^{\infty}(t^{-1}x))=I_{\sigma}^{\infty}(u_{\sigma}^{\infty})=m_{\sigma}^{\infty}. $$ $\hfill\square$

Next, we need the following global compactness lemma of critical type, which is adopted to prove that the functional $I_{\varepsilon,\sigma}$ satisfies $(PS)_{c_{\sigma}}$ condition for a.e. $\sigma\in[\delta,1]$.

{\it\noindent {\bf Lemma 3.7.} If $(V_{2})$, $(F_{3})$, $(F_{4})$ hold and $p\in(2,5).$ For $c>0$ and $\sigma\in[\delta,1],$ and let $ \{u_{n}\}\subset E_{\varepsilon}$ be a bounded $ (PS)_{c}$ sequence for $ I_{\varepsilon,\sigma},$ then there exists $ u\in E_{\varepsilon}$ and $\varrho\in \mathbb R$ such that $ J'_{\varepsilon,\sigma}(u)=0,$ where
$$ J_{\varepsilon,\sigma}(u)=\frac{a+b\varrho^{2}}{2}\int_{\mathbb R^{3}}|\nabla u|^{2}+\frac{1}{2}\int_{\mathbb R^{3}}V(\varepsilon x)|u|^{2}-\frac{\sigma}{p+1}\int_{\mathbb R^{3}}A(\varepsilon x)|u|^{p+1}-\frac{\sigma}{6}\int_{\mathbb R^{3}}B(\varepsilon x)|u|^{6} \eqno(3.23)$$
and either\\
(i) $ u_{n}\rightarrow u$ in $E_{\varepsilon},$ or \\
(ii) there exists $j\in \mathbb N$ and $ \{y_{n}^{k}\}\subset \mathbb R^{3}$ with $|y_{n}^{k}|\rightarrow\infty$ as $ n\rightarrow\infty$ for $ 1\leq k\leq j,$ nontrivial solutions $ w^{1},\cdots,w^{j}$ of the following problem
$$ -(a+b\varrho^{2})\Delta u+V_{\infty}u=\sigma A_{\infty}|u|^{p-1}u+\sigma B_{\infty}|u|^{4}u \eqno(3.24)$$
such that
$$ c+\frac{b\varrho^{4}}{4}=J_{\varepsilon,\sigma}(u)+\sum_{k=1}^{j}J_{\sigma}^{\infty}(w^{k}),$$
where $$ J_{\sigma}^{\infty}(u)=\frac{a+b\varrho^{2}}{2}\int_{\mathbb R^{3}}|\nabla u|^{2}+\frac{1}{2}\int_{\mathbb R^{3}}V_{\infty}|u|^{2}-\frac{\sigma}{p+1}\int_{\mathbb R^{3}}A_{\infty}|u|^{p+1}-\frac{\sigma}{6}\int_{\mathbb R^{3}}B_{\infty}|u|^{6},\eqno(3.25)$$
$$ \left\|u_{n}-u-\sum_{k=1}^{j}w^{k}(\cdot-y_{n}^{k})\right\|_{\varepsilon}\rightarrow0,\ \ \ \ \varrho^{2}=|\nabla u|_{2}^{2}+\sum_{k=1}^{j}|\nabla w^{k}|_{2}^{2}.$$ }

\noindent {\it Proof.} Since $ \{u_{n}\}$ is bounded in $ E_{\varepsilon},$ there exists $ u\in E_{\varepsilon}$ and $ \varrho\in \mathbb R$ such that
$$ u_{n}\rightharpoonup u \ \ in \ \  E_{\varepsilon}, \eqno(3.26)$$
$$ \int_{\mathbb R^{3}}|\nabla u_{n}|^{2}\rightarrow \varrho^{2}. \eqno(3.27)$$
From $(V_{2})$ one deduces immediately that $E_{\varepsilon}$ is embedded continuously into $H^{1}(\mathbb R^{3})$ and therefore also into $L^{q}(\mathbb R^{3}),$ $q\in[2,6]$. Then the fact $I'_{\varepsilon,\sigma}(u_{n})\rightarrow0$ implies that for $\varphi\in C_{0}^{\infty}(\mathbb R^{3}),$
$$ \aligned &\int_{\mathbb R^{3}}(a\nabla u\nabla\varphi+V(\varepsilon x)u\varphi)+b\varrho^{2}\int_{\mathbb R^{3}}\nabla u\nabla \varphi-\sigma\int_{\mathbb R^{3}}A(\varepsilon x)|u|^{p-1}u\varphi-\sigma\int_{\mathbb R^{3}}B(\varepsilon x)|u|^{4}u\varphi=0,\endaligned$$
i.e., $ J'_{\varepsilon,\sigma}(u)=0, $ where
$$J_{\varepsilon,\sigma}(u)=\frac{a+b\varrho^{2}}{2}\int_{\mathbb R^{3}}|\nabla u|^{2}+\frac{1}{2}\int_{\mathbb R^{3}}V(\varepsilon x)|u|^{2}-\frac{\sigma}{p+1}\int_{\mathbb R^{3}}A(\varepsilon x)|u|^{p+1}-\frac{\sigma}{6}\int_{\mathbb R^{3}}B(\varepsilon x)|u|^{6}. $$
Since
$$\aligned  J_{\varepsilon,\sigma}(u_{n})=I_{\varepsilon,\sigma}(u_{n})+\frac{b\varrho^{4}}{4}+o(1),  \endaligned$$
$$ \aligned &\langle J'_{\varepsilon,\sigma}(u_{n}),\varphi \rangle=\langle I'_{\varepsilon,\sigma}(u_{n}),\varphi \rangle+o(1),\ \ \varphi\in E_{\varepsilon},\endaligned  $$
we have
$$J_{\varepsilon,\sigma}(u_{n})\rightarrow c+\frac{b\varrho^{4}}{4}, \ \ \ \ J'_{\varepsilon,\sigma}(u_{n})\rightarrow0  \ \  \hbox{in}\ \ E_{\varepsilon}^{-1}.\eqno(3.28) $$

We claim that either (i) or (ii) holds.

{\bf Step 1.} Set $ u_{n}^{1}=u_{n}-u, $ by (3.26), (3.28), Lemma 2.9, $(V_{2})$, $(F_{3})$ and $(F_{4})$, we have that\\
$ (a.1)$  $ |\nabla u_{n}^{1}|_{2}^{2}=|\nabla u_{n}|_{2}^{2}-|\nabla u|_{2}^{2}+o(1),$ \\
$ (b.1)$  $ |u_{n}^{1}|_{2}^{2}=|u_{n}|_{2}^{2}-|u|_{2}^{2}+o(1),$ \\
$ (c.1)$  $ J_{\sigma}^{\infty}(u_{n}^{1})\rightarrow c+\frac{b\varrho^{4}}{4}-J_{\varepsilon,\sigma}(u),$ \\
$ (d.1)$  $ (J_{\sigma}^{\infty})'(u_{n}^{1})\rightarrow 0 $ in $ E_{\varepsilon}^{-1}.$

Let $$ \sigma^{1}=\limsup_{n\rightarrow\infty}\sup_{y\in \mathbb R^{3}}\int_{B_{1}(y)}(|u_{n}^{1}|^{p+1}+|u_{n}^{1}|^{6}).$$
{\bf Vanishing:}  If $ \sigma^{1}=0,$ applying Lemma 2.8 we conclude that $ u_{n}^{1}\rightarrow0$ in $ L^{s}(\mathbb R^{3}),s\in(2,6]$. By $ (d.1)$, we can obtain that $u_{n}^{1}\rightarrow0$ in $E_{\varepsilon}$ and the proof is completed.\\
{\bf Non-vanishing:} If $\sigma^{1}>0,$ then there exists sequence $ \{y_{n}^{1}\}\subset \mathbb R^{3}$ such that
$$ \int_{B_{1}(y_{n}^{1})}(|u_{n}^{1}|^{p+1}+|u_{n}^{1}|^{6})>\frac{\sigma^{1}}{2}.$$
Set $ w_{n}^{1}\triangleq u_{n}^{1}(\cdot+y_{n}^{1}).$ Then $ \{w_{n}^{1}\}$ is bounded in $ E_{\varepsilon}$ and we may assume that $w_{n}^{1}\rightharpoonup w^{1}$ in $E_{\varepsilon}.$ Hence by (d.1), we obtain $(J_{\sigma}^{\infty})'(w^{1})=0$. Since
$$ \int_{B_{1}(0)}(|w_{n}^{1}|^{p+1}+|w_{n}^{1}|^{6})>\frac{\sigma^{1}}{2},$$
we see that $ w^{1}\neq0.$ Moreover, $u_{n}^{1}\rightharpoonup0$ in $ E_{\varepsilon}$ implies that $\{y_{n}^{1}\}$ is unbounded. Hence, we may assume that $ |y_{n}^{1}|\rightarrow\infty.$

{\bf Step 2.} Set $ u_{n}^{2}=u_{n}-u-w^{1}(\cdot-y_{n}^{1}).$ Similarly, we have\\
$ (a.2)$  $ |\nabla u_{n}^{2}|_{2}^{2}=|\nabla u_{n}|_{2}^{2}-|\nabla u|_{2}^{2}-|\nabla w^{1}|_{2}^{2}+o(1),$ \\
$ (b.2)$  $ |u_{n}^{2}|_{2}^{2}=|u_{n}|_{2}^{2}-|u|_{2}^{2}-|w^{1}|_{2}^{2}+o(1),$ \\
$ (c.2)$  $ J_{\sigma}^{\infty}(u_{n}^{2})\rightarrow c+\frac{b\varrho^{4}}{4}-J_{\varepsilon,\sigma}(u)-J_{\sigma}^{\infty}(w^{1}),$ \\
$ (d.2)$  $ (J_{\sigma}^{\infty})'(u_{n}^{2})\rightarrow 0 $ in $ E_{\varepsilon}^{-1}.$

Let
$$ \sigma^{2}=\limsup_{n\rightarrow\infty}\sup_{y\in \mathbb R^{3}}\int_{B_{1}(y)}(|u_{n}^{2}|^{p+1}+|u_{n}^{2}|^{6}).$$
If vanishing occurs, then $ \|u_{n}^{2}\|_{\varepsilon}\rightarrow0,$ i.e., $ \|u_{n}-u-w^{1}(\cdot-y_{n}^{1})\|_{\varepsilon}\rightarrow0.$ Moreover, by (3.26), (a.2) and (c.2), we have that $$ \varrho^{2}=|\nabla u|_{2}^{2}+|\nabla w^{1}|_{2}^{2}, \ \ \ \ \ c+\frac{b\varrho^{4}}{4}=J_{\varepsilon,\sigma}(u)+J_{\sigma}^{\infty}(w^{1}).$$

If non-vanishing occurs, then there exist sequence $ \{y_{n}^{2}\}\subset \mathbb R^{3}$ and $ w^{2}\in E_{\varepsilon}$ such that $ w_{n}^{2}\triangleq u_{n}^{2}(\cdot+y_{n}^{2})\rightharpoonup w^{2}$ in $ E_{\varepsilon}.$ Then by (d.2), we see that $ (J_{\sigma}^{\infty})'(w^{2})=0.$ Moreover, $ u_{n}^{2}\rightharpoonup0$ in $E_{\varepsilon}$ implies that $|y_{n}^{2}|\rightarrow\infty$ and $ |y_{n}^{2}-y_{n}^{1}|\rightarrow\infty.$

We next proceed by iteration. Recall that if $w^{k}$ is a nontrivial solution of $I_{\sigma}^{\infty},$ then $I_{\sigma}^{\infty}(w^{k})>0.$ Then there exists some finite $j\in \mathbb N$ such that only the vanishing case occurs in Step $ j.$ $\hfill\square$

{\it\noindent {\bf Lemma 3.8.} If $(V_{1})$-$(V_{3})$, $(F_{1})$-$(F_{4})$ hold and $p\in(2,5).$ For $ \sigma\in[\delta,1],$ let $ \{u_{n}\}\subset E_{\varepsilon}$ be a bounded $ (PS)_{c_{\sigma}}$ sequence of $ I_{\varepsilon,\sigma},$ then there exists $ u_{\sigma}\in E_{\varepsilon}$ such that
$ u_{n}\rightarrow u_{\sigma} \ \ in \ \ E_{\varepsilon}.$}

\noindent{\it Proof.} By Lemma 3.7, for $ \sigma\in[\delta,1],$ there exists $ u_{\sigma}\in E_{\varepsilon}$ and $ \varrho_{\sigma}\in \mathbb R$ such that
$$ u_{n}\rightarrow u_{\sigma} \ \ \hbox{in} \ E_{\varepsilon}, \ \ \ \ \int_{\mathbb R^{3}}|\nabla u_{n}|^{2}\rightarrow \varrho_{\sigma}^{2},$$
$ J'_{\varepsilon,\sigma}(u_{\sigma})=0$ and either (i) or (ii) occurs, where $J_{\varepsilon,\sigma}$ is given in (3.23).
If (ii) occurs, i.e., there exists $j\in \mathbb N$ and $ \{y_{n}^{k}\}\subset \mathbb R^{3}$ with $ |y_{n}^{k}|\rightarrow\infty$ as $n\rightarrow\infty$ for $ 1\leq k\leq j,$ nontrivial solutions $ w^{1},\cdots,w^{j}$ of problem (3.24) such that
$$ c_{\sigma}+\frac{b\varrho_{\sigma}^{4}}{4}=J_{\varepsilon,\sigma}(u_{\sigma})+\sum_{k=1}^{j}J_{\sigma}^{\infty}(w^{k}),$$
$$ \left\|u_{n}-u_{\sigma}-\sum_{k=1}^{j}w^{k}(\cdot-y_{n}^{k})\right\|_{\varepsilon}\rightarrow0,\ \ $$
$$ \varrho_{\sigma}^{2}=|\nabla u_{\sigma}|_{2}^{2}+\sum_{k=1}^{j}|\nabla w^{k}|_{2}^{2}, \ \ \ \ \ \ \  \ \ \eqno(3.29)$$
where $J_{\sigma}^{\infty} $ is given in (3.25).

Set
$$\left\{\aligned  &\alpha\triangleq a\int_{\mathbb R^{3}}|\nabla u_{\sigma}|^{2},\  \ \ \ \ \ \ \ \  \ \ \ \ \ \ \ \ \ \ \ \beta\triangleq\int_{\mathbb R^{3}}V(\varepsilon x)|u_{\sigma}|^{2},\\
&\bar{\beta}\triangleq\varepsilon\int_{\mathbb R^{3}}(\nabla V(\varepsilon x),x)|u_{\sigma}|^{2}, \ \ \ \ \  \ \ \ \gamma\triangleq b\varrho_{\sigma}^{2}\int_{\mathbb R^{3}}|\nabla u_{\sigma}|^{2},\\
&\theta\triangleq\int_{\mathbb R^{3}}A(\varepsilon x)|u_{\sigma}|^{p+1}, \ \ \ \ \ \ \ \ \  \ \ \ \ \ \ \bar{\theta}\triangleq\varepsilon\int_{\mathbb R^{3}}(\nabla A(\varepsilon x),x)|u_{\sigma}|^{p+1},\\
&\xi\triangleq\int_{\mathbb R^{3}}B(\varepsilon x)|u_{\sigma}|^{6},\ \ \ \ \ \ \ \ \ \ \ \ \ \  \  \ \ \ \bar{\xi}\triangleq\varepsilon\int_{\mathbb R^{3}}(\nabla B(\varepsilon x),x)|u_{\sigma}|^{6}.\endaligned \right.$$
By Poho\u{z}aev identity applied to $(K_{\varepsilon})$ and $ J'_{\varepsilon,\sigma}(u_{\sigma})=0,$ we have that
$$ \left\{ \aligned &\frac{1}{2}\alpha+\frac{3}{2}\beta+\frac{1}{2}\bar{\beta}+\frac{1}{2}\gamma-\frac{3}{p+1}\sigma\theta-\frac{1}{p+1}\sigma\bar{\theta}-\frac{1}{2}\sigma\xi-\frac{1}{6}\sigma\bar{\xi}=P_{\varepsilon}(u)=0,\\
&\frac{1}{2}\alpha+\frac{1}{2}\beta+\frac{1}{4}\gamma-\frac{\sigma}{p+1}\theta-\frac{1}{6}\sigma\xi=J_{\varepsilon,\sigma}(u_{\sigma})-\frac{1}{4}\gamma,\\
&\alpha+\beta+\gamma-\sigma\theta-\sigma \xi=J'_{\varepsilon,\sigma}(u_{\sigma})=0. \endaligned \right.  $$
Since $\alpha,\gamma,\theta,\xi$ must be nonnegative and $(V_{3}),(F_{1}),(F_{2})$, we obtain
$$\aligned6\left(J_{\varepsilon,\sigma}(u_{\sigma})-\frac{1}{4}\gamma\right)-P_{\varepsilon}(u)-J'_{\varepsilon,\sigma}(u_{\sigma})&=\frac{3}{2}\alpha+\frac{1}{2}(\beta-\bar{\beta})+\frac{p-2}{p+1}\sigma\theta+\frac{1}{p+1}\sigma\bar{\theta}+\frac{1}{2}\sigma\xi+\frac{1}{6}\sigma\bar{\xi}\\
&\geq\frac{p-2}{p+1}\sigma\theta\geq0. \endaligned$$
Hence
$$ J_{\varepsilon,\sigma}(u_{\sigma})\geq\frac{b\varrho_{\sigma}^{2}}{4}\int_{\mathbb R^{3}}|\nabla u_{\sigma}|^{2}. \eqno(3.30)$$
For every nontrivial solution $ w^{k}(k=1,\cdots,j)$ of problem (3.24), i.e., $ (J_{\sigma}^{\infty})'(w^{k})=0.$ Recall that $w^{k}$ satisfies the Poho\u{z}aev identity
$$P_{\sigma}(w^{k})\triangleq\frac{a+b\varrho_{\sigma}^{2}}{2}\int_{\mathbb R^{3}}|\nabla w^{k}|^{2}+\frac{3}{2}\int_{\mathbb R^{3}}V_{\infty}|w^{k}|^{2}-\frac{3\sigma}{p+1}\int_{\mathbb R^{3}}A_{\infty}|w^{k}|^{p+1}-\frac{\sigma}{2}\int_{\mathbb R^{3}}B_{\infty}|w^{k}|^{6}=0, $$
and using (3.29),  we obtain
$$\aligned  0&=\langle(J_{\sigma}^{\infty})'(w^{k}),w^{k}\rangle+P_{\sigma}(w^{k})\\
&=\frac{3(a+b\varrho_{\sigma}^{2})}{2}\int_{\mathbb R^{3}}|\nabla w^{k}|^{2}+\frac{5}{2}\int_{\mathbb R^{3}}V_{\infty}|w^{k}|^{2}-\frac{(p+4)\sigma}{p+1}\int_{\mathbb R^{3}}A_{\infty}|w^{k}|^{p+1}-\frac{3\sigma}{2}\int_{\mathbb R^{3}}B_{\infty}|w^{k}|^{6} \\
&\geq G_{\sigma}^{\infty}(w^{k}).\endaligned $$
Then there exists $ t_{k}\in(0,1]$ such that $ t_{k}w^{k}(t_{k}^{-1}x)\in M_{\sigma}^{\infty}.$ So by (3.29), we have
$$\aligned J_{\sigma}^{\infty}(w^{k})&=\left(J_{\sigma}^{\infty}(w^{k})-\frac{\langle(J_{\sigma}^{\infty})'(w^{k}),w^{k}\rangle
+\bar{Q}(w^{k})}{p+4}-\frac{b\varrho_{\sigma}^{2}}{4}\int_{\mathbb R^{3}}|\nabla w^{k}|^{2}\right)+\frac{b\varrho_{\sigma}^{2}}{4}\int_{\mathbb R^{3}}|\nabla w^{k}|^{2}\\
&\geq I_{\sigma}^{\infty}(t_{k}w^{k}(t_{k}^{-1}x))+\frac{b\varrho_{\sigma}^{2}}{4}\int_{\mathbb R^{3}}|\nabla w^{k}|^{2} \ \ \ \ \ \  \ \ \ \ \ \ \ \ \ \ \ \ \ \ \ \ \ \ \ \ \ \ \ \  \ \ \ \ \ \ \ \ \ \ \ \ \ \ \\
&\geq m_{\sigma}^{\infty}+\frac{b\varrho_{\sigma}^{2}}{4}\int_{\mathbb R^{3}}|\nabla w^{k}|^{2}. \endaligned \eqno(3.31)$$
From (3.29)-(3.31) we see that
$$\aligned  c_{\sigma}+\frac{b\varrho_{\sigma}^{4}}{4}&=J_{\varepsilon,\sigma}(u_{\sigma})+\sum_{k=1}^{j}J_{\sigma}^{\infty}(w^{k})
&\geq m_{\sigma}^{\infty}+\frac{b\varrho_{\sigma}^{4}}{4},
\endaligned$$
i.e., $ c_{\sigma}\geq m_{\sigma}^{\infty},$ which contradicts to Lemma 3.6. So (i) holds, i.e., $ u_{n}\rightarrow u_{\sigma}$ in $ E_{\varepsilon}$ and then $ I_{\varepsilon,\sigma}(u_{\sigma})=c_{\sigma}.$ $\hfill\square$

{\it\noindent {\bf Proof of Theorem 1.1 (existence).}}
We complete the proof in two steps.

{\bf Step 1.}  By Lemma 1.5, Lemma 3.5 and Lemma 3.8, for a.e. $ \sigma\in[\delta,1],$ there exists a nontrivial critical point $ u_{\sigma}\in E_{\varepsilon}$ for $ I_{\varepsilon,\sigma}$ and $I_{\varepsilon,\sigma}(u_{\sigma})=c_{\sigma}.$
Taking a sequence $ \{\sigma_{n}\}\subset[\delta,1]$ satisfying $\sigma_{n}\rightarrow1, $ then we have a sequence nontrivial critical point $\{u_{\sigma_{n}}\}$ of $I_{\varepsilon,\sigma_{n}}$ and $I_{\varepsilon,\sigma_{n}}(u_{\sigma_{n}})=c_{\sigma_{n}}.$ We claim that $\{u_{\sigma_{n}}\}$ is bounded in $ E_{\varepsilon}$. Indeed, let
$$\left\{ \aligned  &\alpha_{n}\triangleq a\int_{\mathbb R^{3}}|\nabla u_{\sigma_{n}}|^{2},\ \ \ \ \ \ \ \ \ \ \ \ \ \ \ \ \ \ \beta_{n}\triangleq\int_{\mathbb R^{3}}V(\varepsilon x)|u_{\sigma_{n}}|^{2},\\
&\bar{\beta}_{n}\triangleq\varepsilon\int_{\mathbb R^{3}}(\nabla V(\varepsilon x),x)|u_{\sigma_{n}}|^{2},\ \ \ \ \ \ \gamma_{n}\triangleq b\left(\int_{\mathbb R^{3}}|\nabla u_{\sigma_{n}}|^{2}\right)^{2},\\
&\theta_{n}\triangleq\int_{\mathbb R^{3}}A(\varepsilon x)|u_{\sigma_{n}}|^{p+1},\ \ \ \ \ \ \ \ \ \ \ \ \ \bar{\theta}_{n}\triangleq\varepsilon\int_{\mathbb R^{3}}(\nabla A(\varepsilon x),x)|u_{\sigma_{n}}|^{p+1},\\
&\xi_{n}\triangleq\int_{\mathbb R^{3}}B(\varepsilon x)|u_{\sigma_{n}}|^{6},\ \ \ \ \ \ \ \ \ \ \ \ \ \ \ \ \bar{\xi}_{n}\triangleq\varepsilon\int_{\mathbb R^{3}}(\nabla B(\varepsilon x),x)|u_{\sigma_{n}}|^{6},\endaligned \right.$$
we get
$$\left\{ \aligned &\frac{1}{2}\alpha_{n}+\frac{3}{2}\beta_{n}+\frac{1}{2}\bar{\beta_{n}}+\frac{1}{2}\gamma_{n}-\frac{3}{p+1}\sigma_{n}\theta_{n}-\frac{1}{p+1}\sigma_{n}\bar{\theta}_{n}-\frac{1}{2}\sigma_{n}\xi_{n}-\frac{1}{6}\sigma_{n}\bar{\xi}_{n}=P_{\varepsilon}(u)=0,\\
&\frac{1}{2}\alpha_{n}+\frac{1}{2}\beta_{n}+\frac{1}{4}\gamma_{n}-\frac{1}{p+1}\sigma_{n}\theta_{n}-\frac{1}{6}\sigma_{n}\xi_{n}=I_{\varepsilon,\sigma_{n}}(u_{\sigma_{n}})=c_{\sigma_{n}},\\
&\alpha_{n}+\beta_{n}+\gamma_{n}-\sigma_{n}\theta_{n}-\sigma_{n}\xi_{n}=I'_{\varepsilon,\sigma_{n}}(u_{\sigma_{n}})=0. \endaligned \right.  $$
So that
$$\aligned 6c_{\delta}\geq 6c_{\sigma_{n}}&=6c_{\sigma_{n}}-P_{\varepsilon}(u)-I'_{\varepsilon,\sigma_{n}}(u_{\sigma_{n}})\\
&=\frac{3}{2}\alpha_{n}+\frac{1}{2}\left(\beta_{n}-\bar{\beta}_{n}\right)+\frac{p-2}{p+1}\sigma_{n}\theta_{n}+\frac{1}{p+1}\sigma_{n}\bar{\theta}_{n}+\frac{1}{2}\sigma_{n}\xi_{n}+\frac{1}{6}\sigma_{n}\bar{\xi}_{n}, \endaligned\eqno(3.32) $$
$$c_{\sigma_{n}}=c_{\sigma_{n}}-\frac{1}{4}I'_{\varepsilon,\sigma_{n}}(u_{\sigma_{n}})=\frac{1}{4}(\alpha_{n}+\beta_{n})+\left(\frac{1}{4}-\frac{1}{p+1}\right)\sigma_{n}\theta_{n}+\frac{1}{12}\sigma_{n}\xi_{n}. \eqno(3.33) $$
We then obtain immediately, by the assumptions $(V_{3}), (F_{1}), (F_{2})$, $\beta_{n}-\bar{\beta}_{n}\geq0, \bar{\theta}_{n}\geq0, \bar{\xi}_{n}\geq0$ and $ \alpha_{n}, \beta_{n}, \theta_{n}, \xi_{n}>0$. This means that $ \theta_{n}, \xi_{n}$ are bounded from (3.32), hence by (3.33), $\alpha_{n}+\beta_{n} $ is bounded, i.e., $\{ u_{\sigma_{n}}\}$ is bounded in $E_{\varepsilon}.$ Then by Lemma 3.5, we see that
$$ \lim_{n\rightarrow\infty}I_{\varepsilon,1}(u_{\sigma_{n}})=\lim_{n\rightarrow\infty}\left(I_{\varepsilon,\sigma_{n}}(u_{\sigma_{n}})+(\sigma_{n}-1)\int_{\mathbb R^{3}}
|u_{\sigma_{n}}|^{p+1}\right)=\lim_{n\rightarrow\infty}c_{\sigma_{n}}=c_{1},$$
$$ \lim_{n\rightarrow\infty}\langle I'_{\varepsilon,1}(u_{\sigma_{n}}),\varphi\rangle=\lim_{n\rightarrow\infty}\left(\langle I'_{\varepsilon,\sigma_{n}}(u_{\sigma_{n}}),\varphi\rangle+(\sigma_{n}-1)\int_{\mathbb R^{3}}
|u_{\sigma_{n}}|^{p-1}u_{\sigma_{n}}\varphi \right)=0,\ \forall \varphi\in E_{\varepsilon},$$
i.e., $ \{u_{\sigma_{n}}\}$ is a bounded $(PS)_{c_{1}}$ sequence for $ I_{\varepsilon}=I_{\varepsilon,1}.$ Therefore by Lemma 3.8, there exists a nontrivial critical point $ u_{1}\in E_{\varepsilon}$ for $I_{\varepsilon} $ and $I_{\varepsilon}(u_{1})=c_{1}.$

{\bf Step 2.} Next we prove the existence of a ground state solution of Nehari-Poho\u{z}aev type for problem $(1.1)$.
Set
$$ m=\inf_{X}I_{\varepsilon}(u),\ \ \hbox{where}\ X:=\{u\in E_{\varepsilon}\backslash\{0\}, I'_{\varepsilon}(u)=0\}.$$
Then we see that $0<m\leq I_{\varepsilon}(u_{1})=c_{1}<+\infty.$  Indeed, for all $u\in X$, we have $I'_{\varepsilon}(u)u=0$, i.e.,
$$\|u\|_{\varepsilon}^{2}+b\left(\int_{\mathbb R^{3}}|\nabla u|^{2}\right)^{2}-\int_{\mathbb R^{3}}A(\varepsilon x)|u|^{p+1}-\int_{\mathbb R^{3}}B(\varepsilon x)|u|^{6}=0.$$
and
$$\aligned\|u\|_{\varepsilon}^{2}&\leq\int_{\mathbb R^{3}}A(\varepsilon x)|u|^{p+1}+\int_{\mathbb R^{3}}B(\varepsilon x)|u|^{6}&\leq A_{0}|u|_{p+1}^{p+1}+B_{0}|u|_{6}^{6}\leq 2C\|u\|_{\varepsilon}^{6}.\endaligned$$
So for some $C> 0$, $\|u\|_{\varepsilon}\geq C>0.$ On the other hand, the Poho\u{z}aev identity (2.1) holds, i.e., $P_{\varepsilon}(u)=0.$
Now we can calculate
$$\aligned I_{\varepsilon}(u)&=I_{\varepsilon}(u)-\frac{1}{3}P_{\varepsilon}(u)\\
&=\frac{a}{3}\int_{\mathbb R^{3}}|\nabla u|^{2}-\frac{\varepsilon}{6}\int_{\mathbb R^{3}}(\nabla V(\varepsilon x),x)u^{2}+\frac{1}{12}b\left(\int_{\mathbb R^{3}}|\nabla u|^{2}\right)^{2}\\
&\ \ \ +\frac{\varepsilon}{3(p+1)}\int_{\mathbb R^{3}}(\nabla A(\varepsilon x),x)|u|^{p+1}+\frac{\varepsilon}{6}\int_{\mathbb R^{3}}(\nabla B(\varepsilon x),x)|u|^{6}. \endaligned$$
By $(V_{4})$ and Hardy inequality
$$ \int_{\mathbb R^{3}}|\nabla u|^{2}\geq \frac{1}{4}\int_{\mathbb R^{3}}\frac{|u|^{2}}{|x|^{2}},$$
we infer
$$I_{\varepsilon}(u)\geq \frac{a-l}{12}\int_{\mathbb R^{3}}|\nabla u|^{2}\geq0.\eqno(3.34)$$
Therefore, we obtain $m\geq0.$ In the following let us rule out $m=0.$ By contradiction, let $\{u_{n}\}$ be a $(PS)_{0}$ sequence of $I_{\varepsilon}$. From (3.34), we have $\lim\limits_{n\rightarrow\infty}\int_{\mathbb R^{3}}|\nabla u_{n}|^{2}=0.$ This conclusion combined with $I'_{\varepsilon}(u)u=0$ imply $\lim\limits_{n\rightarrow\infty}\int_{\mathbb R^{3}}V(\varepsilon x)|u_{n}|^{2}=0.$ Therefore, we obtain $\lim\limits_{n\rightarrow\infty}\|u_{n}\|_{\varepsilon}=0,$ a contradiction with $\|u\|_{\varepsilon}\geq C>0$ for all $n\in \mathbb N.$
Let $ \{u_{n}\}$ be a sequence of nontrivial critical points of $ I_{\varepsilon}$ satisfying $ I_{\varepsilon}(u_{n})\rightarrow m,$ using the same arguments as in {\bf Step 1}, we conclude that $\{u_{n}\}$ is bounded in $E_{\varepsilon},$ i.e., $ \{u_{n}\}$ is a bounded $ (PS)_{m}$ sequence of $ I_{\varepsilon}.$ By Lemma 3.8 and strong maximum principle, there exists a nontrivial $ u\in E_{\varepsilon}$ such that $I_{\varepsilon}(u)=m$ and
$I'_{\varepsilon}(u)=0$ and $u(x)>0$ for all $x\in \mathbb R^{3}.$ So $u$ is a positive ground state solution for problem $(1.1)$. The proof is completed. $\hfill\square$

\vskip4pt

\section{ Concentration of positive ground state of Nehari-Poho\u{z}aev type}

\indent In this section, we will adopt the family $v_{\varepsilon}(x)=u_{\varepsilon}(\varepsilon x)$ to study the behavior of the solution $u_{\varepsilon}(x)$ of $(SK_{\varepsilon}).$ Without loss of generality, we assume $0\in\Lambda\cap \Lambda_{1}\cap \Lambda_{2}$. Here, we give the associated autonomous problem
$$\left\{\aligned &-\left(a+b\int_{\mathbb{R}^{3}}|\nabla u|^{2}\right)\Delta u+V_{0}u=A_{0}|u|^{p-1}u+B_{0}|u|^{4}u, &x \in \mathbb{R}^{3},\\
&u\in H^{1}(\mathbb{R}^{3}),\ \  u>0, &x \in \mathbb{R}^{3},\endaligned\right.  \eqno(K_{0})$$
which will be used to prove the concentration of positive ground state solutions. We give the following energy functional associated to $(K_{0})$
$$ I_{0}(u)=\frac{a}{2}\int_{\mathbb R^{3}}|\nabla u|^{2}+\frac{1}{2}\int_{\mathbb R^{3}}V_{0}|u|^{2}+\frac{b}{4}\left(\int_{\mathbb R^{3}}|\nabla u|^{2}\right)^{2}-\frac{1}{p+1}\int_{\mathbb R^{3}}A_{0}|u|^{p+1}-\frac{1}{6}\int_{\mathbb R^{3}}B_{0}|u|^{6},\eqno(4.1)$$
and set
$$ M_{0}=\{u\in H^{1}(\mathbb R^{3})\backslash\{0\}:G_{0}(u)=0\},\ \ \  m_{0}=\inf\limits_{u\in M_{0}}I_{0}(u),$$
where $G_{0}(u)=\langle I'_{0}(u),u \rangle+P_{0}(u)$,
$$P_{0}(u):=\frac{1}{2}a\int_{\mathbb R^{3}}|\nabla u|^{2}+\frac{3}{2}\int_{\mathbb R^{3}}V_{0}|u|^{2}+\frac{1}{2}b\left(\int_{\mathbb R^{3}}|\nabla u|^{2}\right)^{2}-\frac{3}{p+1}\int_{\mathbb R^{3}}A_{0}|u|^{p+1}-\frac{1}{2}\int_{\mathbb R^{3}}B_{0}|u|^{6}=0$$ is the Poho\u{z}aev equality applied to (4.1).

{\it\noindent {\bf Lemma 4.1.}  $(K_{0})$ has a positive ground state solution of Nehari-Poho\u{z}aev type in $H^{1}(\mathbb R^{3})$. In particular, there exists a minimizer of $ m_{0}$.}

{\it\noindent Proof.} According to Section 2 and Section 3, $M_{0}$ has some similar properties to those of the manifold $M$, such as containing all the nontrivial critical points of $I_{0}$. The proof which is similar to Theorem 1.6, and is omitted here. $\hfill\square$

{\it\noindent {\bf Lemma 4.2.} $\lim\limits_{\varepsilon\rightarrow0^{+}}c_{\varepsilon}=m_{0}$.}\\
\noindent {\it Proof.} Let $\kappa\in C^{\infty}(\mathbb R^{3},[0,1])$ be a function such that $\kappa(x)=1$ if $|x|\leq\frac{1}{2}$ and $\kappa(x)=0$ if $|x|\geq1.$ For every $R>0$ and every $x\in \mathbb R^{3}$, set $u_{R}(x)=\Psi_{R}(x)u_{0},$ here $u_{0}$ is a positive ground state solution of Nehari-Poho\u{z}aev type of $(K_{0})$ and $\Psi_{R}(x)=\kappa(\frac{x}{R})$.  By Lebesgue Theorem, we have
$$ u_{R}\rightarrow u_{0}\ \ \hbox{as}\ R\rightarrow\infty. \eqno(4.2)$$
For every $\varepsilon,R>0$, there exists $t_{\varepsilon,R} > 0$ such that
$$I_{\varepsilon}(t_{\varepsilon,R}u_{R})=\max_{t\geq0}I_{\varepsilon}(tu_{R}).$$
Hence,
$$\aligned &\ \frac{1}{t_{\varepsilon,R}^{2}}\int_{B_{R}(0)}\left(a|\nabla u_{R}|^{2}+V(\varepsilon x\right)|u_{R}|^{2})+b\left( \int_{B_{R}(0)}|\nabla u_{R}|^{2}\right)^{2}\\
=&\ t_{\varepsilon,R}^{p-3}\int_{B_{R}(0)}A(\varepsilon x)|u_{R}|^{p+1}+t_{\varepsilon,R}^{2}\int_{B_{R}(0)}B(\varepsilon x)|u_{R}|^{6},\ \ \ \ \ \ \  \ \ \ \ \ \  \endaligned \eqno(4.3)$$
which implies that
$$\aligned &\frac{1}{t_{\varepsilon,R}^{2}}\int_{B_{R}(0)}(a|\nabla u_{R}|^{2}+|V|_{\infty(|x|<R)}|u_{R}|^{2})+b\left( \int_{B_{R}(0)}|\nabla u_{R}|^{2}\right)^{2}\\
\geq &\ t_{\varepsilon,R}^{p-3}\int_{B_{R}(0)}A_{\infty}|u_{R}|^{p+1}+t_{\varepsilon,R}^{2}\int_{B_{R}(0)}B_{\infty}|u_{R}|^{6}. \endaligned$$
From (4.3) and the last inequality we see that
$$ 0<\lim_{\varepsilon\rightarrow0^{+}}t_{\varepsilon,R}=t_{R}<\infty,$$
for every $R>0$.
Let $\varepsilon\rightarrow0^{+}$ in (4.3), we get
$$\aligned &\frac{1}{t_{R}^{2}}\int_{B_{R}(0)}\left(a|\nabla u_{R}|^{2}+V_{0}|u_{R}|^{2}\right)+b\left(\int_{B_{R}(0)}|\nabla u_{R}|^{2}\right)^{2}\\
=&\ t_{R}^{p-3}\int_{B_{R}(0)}A_{0}|u_{R}|^{p+1}+t_{R}^{2}\int_{B_{R}(0)}B_{0}|u_{R}|^{6}. \endaligned \eqno(4.4)$$
Combing (4.2) with (4.4), we get $\lim \limits_{R\rightarrow\infty}t_{R}=1$ and $I_{0}(t_{R}u_{R})=\max\limits_{t\geq0}I_{0}(tu_{R})$.
Therefore, $c_{\varepsilon}\leq \max\limits_{t\geq0}I_{\varepsilon}(tu_{R})=I_{\varepsilon}(t_{\varepsilon,R}u_{R})$ implies
$ \limsup\limits_{\varepsilon\rightarrow0^{+}}c_{\varepsilon}\leq I_{0}(t_{R}u_{R}).$ Using (4.2), we deduce that $ \limsup\limits_{\varepsilon\rightarrow0^{+}}c_{\varepsilon}\leq m_{0}.$ In order to prove $\liminf\limits_{\varepsilon\rightarrow0^{+}}c_{\varepsilon}\geq m_{0}$, it is suffices to verify that $c_{\varepsilon}\geq m_{0}$ for all $\varepsilon\in (0, \varepsilon^{*})$. In fact, we assume on the contrary that $c_{\varepsilon_{0}^{*}}<m_{0}$ for some $\varepsilon_{0}^{*}\in (0, \varepsilon^{*}).$ From the part (1) of Theorem 1.1, there exists $u_{\varepsilon_{0}^{*}}\in E_{\varepsilon_{0}^{*}}\backslash\{0\}$ such that $c_{\varepsilon_{0}^{*}}=I_{\varepsilon_{0}^{*}}(u_{\varepsilon_{0}^{*}})=\max\limits_{t>0}I_{\varepsilon_{0}^{*}}(tu_{\varepsilon_{0}^{*}})<m_{0}.$ By the definition of $m_{0}$, we see that $m_{0}\leq \max\limits_{t>0}I_{0}(tu_{\varepsilon_{0}^{*}})$. Since $V_{0}\leq V(\varepsilon_{0}^{*}x)$, $A_{0}\geq A(\varepsilon_{0}^{*}x)$, $B_{0}\geq B(\varepsilon_{0}^{*}x)$, we have $m_{0}>\max\limits_{t>0}I_{\varepsilon_{0}^{*}}(tu_{\varepsilon_{0}^{*}})\geq \max\limits_{t>0}I_{0}(tu_{\varepsilon_{0}^{*}})\geq m_{0},$ which is impossible and therefore $\lim\limits_{\varepsilon\rightarrow0^{+}}c_{\varepsilon}=m_{0}.$ $\hfill\square$

{\it\noindent {\bf Lemma 4.3.} Let $v_{n}(x)=v_{\varepsilon_{n}}(x)$ and $\{v_{n}\}\subset M_{0}$  be a sequence satisfying $ I_{0}(v_{n})\rightarrow m_{0}$ as $ n\rightarrow\infty.$ If $(V_{2})$ and $(F_{5})$ hold, for the family $ v_{\varepsilon}$ satisfying $I_{\varepsilon}(v_{\varepsilon})=c_{\varepsilon}$ and $I'_{\varepsilon}(v_{\varepsilon})=0$, there exist $ \varepsilon^{*}>0$, a family $ \{y_{\varepsilon}\}\subset \mathbb R^{3}$  such that the family $ \{\varepsilon y_{\varepsilon}\}$ is bounded. In particular, if $ x_{0}$ is limit of the sequence $\{\varepsilon_{n} y_{\varepsilon_{n}}\}$ in the family $\{\varepsilon y_{\varepsilon}\}, $ then we have $ V(x_{0})=V_{0}$ and $x_{0}\in \Lambda \cap \Lambda_{1} \cap \Lambda_{2}.$}\\
\noindent {\it Proof.}  Using Lemma 4.1, $w_{n}\rightarrow w$ in $H^{1}(\mathbb R^{3})$. By Fatou's Lemma, $G_{0}(w)=0$ and $ \{w_{n}\}\subset M_{0}$, we know that
$$\aligned m_{0}&\leq I_{0}(w)<I_{\infty}(w)-\frac{1}{6}G_{0}(w)\\
&=\frac{a}{4}\int_{\mathbb R^{3}}|\nabla w|^{2}+\frac{1}{12}\int_{\mathbb R^{3}}(V_{\infty}-V_{0})|w|^{2}+\frac{p-2}{6(p+1)}\int_{\mathbb R^{3}}(A_{0}-A_{\infty})|w|^{p+1}\\
&\ \ \ +\frac{1}{12}\int_{\mathbb R^{3}}\left(B_{0}-B_{\infty}\right)|w|^{6}\\
&\leq\liminf_{n\rightarrow\infty}\bigg(\frac{a}{4}\int_{\mathbb R^{3}}|\nabla w_{n}|^{2}+\frac{1}{12}\int_{\mathbb R^{3}}(V_{\infty}-V_{0})|w_{n}|^{2}+\frac{p-2}{6(p+1)}\int_{\mathbb R^{3}}(A_{0}-A_{\infty})|w_{n}|^{p+1}\\
&\ \ \ +\frac{1}{12}\int_{\mathbb R^{3}}\left(B_{0}-B_{\infty}\right)|w_{n}|^{6}\bigg)\\
&=\liminf_{n\rightarrow\infty}I_{\varepsilon_{n}}(v_{n})\\
&=m_{0}, \endaligned \eqno(4.5)$$
which is a contradiction. Hence, $\{\varepsilon_{n}y_{n} \} $ is bounded and up to a subsequence, $\varepsilon_{n}y_{n}\rightarrow x_{0}$.
Consider the functional $I_{x_{0}}: H^{1}(\mathbb R^{3})\rightarrow\mathbb R$ as
$$I_{x_{0}}(u):=\frac{1}{2}\int_{\mathbb R^{3}}(a|\nabla u|^{2}+V(x_{0})u^{2})+\frac{b}{4}\left(\int_{\mathbb R^{3}}|\nabla u|^{2}\right)^{2}-\int_{\mathbb R^{3}}A(x_{0})|u|^{p+1}+\frac{1}{6}\int_{\mathbb R^{3}}B(x_{0})|u|^{6}.$$
Then, if $V(x_{0})>V_{0}$, we can get a contradiction by using $I_{x_{0}}$ to take the place of $ I_{\infty}$ in (4.5). It follows from $(V_{2})$ that $V(x_{0})=V_{0}$, i.e., $ x_{0}\in \Lambda$. By $(F_{5})$, if $x_{0}\not\in \Lambda_{1}\cap\Lambda_{2}$, that is, $ A(x_{0})<A_{0}$ or $B(x_{0})<B_{0}$, then $ I_{0}(\bar{w})<I_{x_{0}}(\bar{w}).$ Hence, we can also obtain a contradiction by repeating the same arguments used above. Therefore, $ x_{0}\in \Lambda\cap\Lambda_{1}\cap\Lambda_{2}$. $\hfill\square$

In the following, we introduce two Lemmas, which is a key point in the study of concentration of positive ground state solutions of Nehari-Poho\u{z}aev type of $(SK_{\varepsilon})$.

\noindent {\bf Lemma 4.4.} (\cite{40}) {\it Let $ u\in H^{1}(\mathbb R^{3})$ satisfies\\
\indent \indent  \indent $ -\Delta u+(b(x)-q(x))u=h(x,u)$  in  $ \mathbb R^{3}$,\\
and $ h:\mathbb R^{3}\times \mathbb R\rightarrow \mathbb R^{+}$ be a Caratheodory function satisfying\\
\indent \indent \indent $ 0\leq h(x,v)\leq C_{h}(v^{r}+v)$ \ for all $v>0, \ x\in \mathbb R^{3},$ \\
where $b:\mathbb R^{3}\rightarrow \mathbb R^{+}$ is a $L_{loc}^{\infty}(\mathbb R^{3})$ function, $q(x)\in L^{\frac{3}{2}}(\mathbb R^{3}), \ 1<r<5.$ Then, $u\in L^{s}(\mathbb R^{3})$ for all $s\geq2.$ Moreover, there is a positive constant $ C_{s}$ depending on $s,$ $q(x)$ and $C_{h},$ such that $|u|_{s}\leq C_{s}\|u\|_{H}.$}

\noindent {\bf Lemma 4.5.} (\cite{41}) {\it Suppose that $t>3, \ g\in L^{\frac{t}{2}}(\Omega)$ and $u\in H^{1}(\Omega)$ satisfy in the weak sense\\
\indent \indent \indent \indent \indent \indent \indent \indent$-\Delta u\leq g(x)$ \  in $\Omega,$ \\
where $\Omega$ is an open subset of $\mathbb R^{3}.$ Then, for every ball $B_{2R}(y)\subset\Omega, R>0,$
$$ \sup_{x\in B_{R}(y)}u(x)\leq C(|u^{+}|_{L^{2}(B_{R}(y))}+|g|_{L^{2}(B_{R}(y))}),$$
where $C$ depends on $t$ and $R.$}

{\it\noindent {\bf Lemma 4.6.} $w_{n}\rightarrow w$ in $H^{1}(\mathbb R^{3})$ as $n\rightarrow\infty,$ where the sequence $\{w_{n}\}$ has been given in Lemma 4.3.
Furthermore, there exists $\varepsilon^{*}>0$ such that $\lim\limits_{|x|\rightarrow\infty}w_{\varepsilon}(x)=0$ uniformly on $\varepsilon\in(0,\varepsilon^{*}).$}\\
\noindent {\it Proof.} It follows from Lemma 4.3 that $w_{n}\rightarrow w$ in $H^{1}(\mathbb R^{3})$. Hence, we have
$$ \lim_{R\rightarrow\infty}\int_{|x|\geq R}(w_{n}^{2}+w_{n}^{6})=0,\ \ \hbox{for} \ n\in \mathbb N. \eqno(4.6)$$
Let $B^{1}_{n}=a+b\int_{\mathbb R^{3}}|\nabla w_{n}|^{2},$ then there exists $c>0$ such that $0<a<B^{1}_{n}\leq c$. It follows from Lemma 4.3 that $w_{n}$ solves
$$ -\triangle w_{n}+\frac{V(\varepsilon_{n}x+\varepsilon_{n}y_{n})}{B^{1}_{n}}w_{n}=
\frac{A(\varepsilon_{n}x+\varepsilon_{n}y_{n})}{B^{1}_{n}}|w_{n}|^{p-1}w_{n}+\frac{B(\varepsilon_{n}x+\varepsilon_{n}y_{n})}{B^{1}_{n}}|w_{n}|^{4}w_{n}.$$
By Lemma 4.4 with $q(x)=\frac{B(\varepsilon_{n}x+\varepsilon_{n}y_{n})}{B^{1}_{n}}|w_{n}|^{4}\in L^{\frac{3}{2}}(\mathbb R^{3}),\ b(x)=\frac{V(\varepsilon_{n}x+\varepsilon_{n}y_{n})}{B^{1}_{n}},\ h(x,w_{n})=\frac{A(\varepsilon_{n}x+\varepsilon_{n}y_{n})}{B^{1}_{n}}|w_{n}|^{p-1}w_{n},$ we obtain $ w_{n}\in L^{t}(\mathbb R^{3})$ for all $t\geq2$ and $|w_{n}|_{t}\leq C_{t}\|w_{n}\|_{H}.$
Set
$$g_{n}(x)=\frac{A(\varepsilon_{n}x+\varepsilon_{n}y_{n})}{B^{1}_{n}}|w_{n}|^{p-1}w_{n}+\frac{B(\varepsilon_{n}x+\varepsilon_{n}y_{n})}{B^{1}_{n}}|w_{n}|^{4}w_{n}, $$
then, we have $|g_{n}|_{\frac{t}{2}}\leq C $ for $t>3$ and all $n.$ Using Lemma 4.5, we have
$$ \sup_{x\in B_{1}(y)}w_{n}(x)\leq C(|w_{n}|_{L^{2}(B_{1}(y))}+|g_{n}|_{L^{2}(B_{1}(y))}),\ \ \hbox{ for  all } \ y\in \mathbb R^{3}, \eqno(4.7)$$
which implies that $|w_{n}|_{\infty}$ is uniformly bounded for $n\in \mathbb N$. Combing (4.6) with (4.7), we see that $ \lim\limits_{|x|\rightarrow\infty}w_{n}(x)=0$ uniformly  on $n\in \mathbb N.$
Consequently, there exists $\varepsilon^{*}>0$ such that
$$ \lim_{|x|\rightarrow\infty}w_{\varepsilon}(x)=0\ \ \hbox{uniformly \ on} \ \varepsilon\in(0,\varepsilon^{*}).$$  $\hfill\square$

\indent Now we consider the exponential decay of solutions $w_{\varepsilon}$.

{\it\noindent {\bf Lemma 4.7.}  There exist constants $C>0$ and $c>0 $ such that
$$ w_{\varepsilon}(x)\leq Ce^{-c|x|}, \ \ for \ all \ x\in \mathbb R^{3}. \eqno(4.8)$$}
\noindent {\it Proof.} By Lemma 4.6, there exist $R,\varepsilon^{*}>0$ such that
$$ \frac{A(\varepsilon x+\varepsilon y_{\varepsilon})}{(a+b\int_{\mathbb R^{3}}|\nabla w_{\varepsilon}|^{2})}|w_{\varepsilon}|^{p-1}+\frac{B(\varepsilon x+\varepsilon y_{\varepsilon})}{(a+b\int_{\mathbb R^{3}}|\nabla w_{\varepsilon}|^{2})}|w_{\varepsilon}|^{4}\leq \frac{V_{0}}{2(a+bM)}$$
for all $|x|\geq R$ and $\varepsilon\in(0,\varepsilon^{*})$, where $ M\geq \|w_{\varepsilon}\|_{H}^{2}$  uniformly for $\varepsilon\in(0,\varepsilon^{*})$.
Fix $ \varphi(x)=Ce^{-c|x|}$ with $ c^{2}<\frac{V_{0}}{2(a+bM)}$ and $Ce^{-c|x|}\geq w_{\varepsilon}$ for all $|x|=R$. Then
$$ \Delta\varphi\leq c^{2}\varphi, \ \ \hbox{for all } \ x\neq0. \eqno(4.9)$$
Since $w_{\varepsilon}>0$, we have
$$ \aligned-\Delta w_{\varepsilon}+\frac{V_{0}}{a+bM}w_{\varepsilon}&\leq-\Delta w_{\varepsilon}+\frac{V_{0}}{(a+b\int_{\mathbb R^{3}}|\nabla w_{\varepsilon}|^{2})}w_{\varepsilon}\\
&\leq\frac{A(\varepsilon x+\varepsilon y_{\varepsilon})}{(a+b\int_{\mathbb R^{3}}|\nabla w_{\varepsilon}|^{2})}|w_{\varepsilon}|^{p-1}w_{\varepsilon}+\frac{B(\varepsilon x+\varepsilon y_{\varepsilon})}{(a+b\int_{\mathbb R^{3}}|\nabla w_{\varepsilon}|^{2})}|w_{\varepsilon}|^{4}w_{\varepsilon}\\
&\leq\frac{V_{0}}{2(a+bM)}w_{\varepsilon}, \ \ \hbox{ for } |x|\geq R.\endaligned \eqno(4.10)$$
Let $ \varphi_{\varepsilon}=\varphi-w_{\varepsilon}$. Combining (4.9) with (4.10), we have
$$ \left\{\aligned &-\Delta\varphi_{\varepsilon}+\frac{V_{0}}{2(a+bM)}\varphi_{\varepsilon}\geq0 \  \ \hbox{in} \ |x|\geq R,\\
& \varphi_{\varepsilon}\geq0 \ \ \hbox{on} \ |x|=R,\\
&\lim_{|x|\rightarrow\infty}\varphi_{\varepsilon}(x)=0. \endaligned \right. $$
The maximum principle implies that $\varphi_{\varepsilon}\geq0$ in $|x|\geq R$, and we obtain
$$w_{\varepsilon}(x)\leq Ce^{-c|x|} \ \hbox{ for all } |x|\geq R \ \hbox{and} \ \varepsilon\in(0,\varepsilon^{*}).$$ $\hfill\square$

{\it\noindent {\bf Proof of Theorem 1.1 (concentration).}} From Lemma 4.6 and $w_{\varepsilon}=v_{\varepsilon}(x+y_{\varepsilon}),$ there exists $R>0$ such that the global maximum point of $w_{\varepsilon},$ denoted by $k_{\varepsilon},$ is located in $B_{R}(0)$. Then, the global maximum point of $v_{\varepsilon}$, given by $z_{\varepsilon}=y_{\varepsilon}+k_{\varepsilon},$ satisfies $\varepsilon z_{\varepsilon}=\varepsilon y_{\varepsilon}+\varepsilon k_{\varepsilon}.$ Notice that $u_{\varepsilon}(x)=v_{\varepsilon}(\frac{x}{\varepsilon}),$ then we see that $u_{\varepsilon}(x)$ is a positive ground state solution of Nehari-Poho\u{z}aev type to $(SK_{\varepsilon})$ with $\varepsilon>0$ and has a global maximum point $x_{\varepsilon}=\varepsilon z_{\varepsilon}$. By Lemma 4.3 and Lemma 4.6, we have $ \lim\limits_{\varepsilon\rightarrow0}V(x_{\varepsilon})=V_{0}.$
From the definition of $u_{\varepsilon}$ and (4.8), we deduce that
$$ u_{\varepsilon}(x)=v_{\varepsilon}\left(\frac{x}{\varepsilon}\right)=w_{\varepsilon}\bigg(\frac{x}{\varepsilon}-y_{\varepsilon}\bigg)=w_{\varepsilon}\bigg(\frac{x+\varepsilon k_{\varepsilon}-x_{\varepsilon}}{\varepsilon}\bigg)\leq C\exp\left(-c\frac{|x-x_{\varepsilon}|}{\varepsilon}\right)$$
for all $x\in \mathbb R^{3}$ and $ \varepsilon\in(0,\varepsilon^{*})$.  $\hfill\square$

\noindent{\bf Availability of data and materials}

All data generated or analyzed during this study are included in this paper.

\noindent{\bf Acknowledgments}

The authors thank anonymous referees whose important comments helped them to improve their work.
The  authors were supported by NSFC (11571197).

\end{CJK*}
\end{document}